\documentclass[12]{amsart}
\usepackage{amsmath, amssymb, verbatim, graphicx, amsthm, color}

\newtheorem{theorem}{Theorem}[section]
\newtheorem{corollary}{Corollary}[section]
\newtheorem{proposition}{Proposition}[section]
\newtheorem{lemma}{Lemma}[section]

\newtheorem*{ebtheorem}{Theorem}

\numberwithin{equation}{section}

\newcommand{\CH}{\mathcal H}
\newcommand{\CL}{\mathcal L}
\newcommand{\CB}{\mathcal B}

\newcommand{\CN}{\mathcal N}
\newcommand{\CP}{\mathcal P}
\newcommand{\CQ}{\mathcal Q}

\newcommand{\CU}{\mathcal U}

\newcommand{\R}{\mathbb R}
\newcommand{\Z}{\mathbb Z}
\newcommand{\T}{\mathbb T}
\newcommand{\TT}{\mathbb T}

\newcommand{\trace}{{\rm trace}}
\newcommand{\grad}{{\rm grad}}
\newcommand{\Div}{{\rm div}}
\newcommand{\diff}{{\rm Diff}}

\newcommand{\Endo}{\mathrm{Endo}}

\newcommand{\beq}{\begin{equation}}
\newcommand{\eeq}{\end{equation}}
\newcommand{\beqn}{\begin{equation*}}
\newcommand{\eeqn}{\end{equation*}}
\newcommand{\beqy}{\begin{eqnarray}}
\newcommand{\eeqy}{\end{eqnarray}}
\newcommand{\beqyn}{\begin{eqnarray*}}
\newcommand{\eeqyn}{\end{eqnarray*}}

\newcommand{\p}{\partial}

\newcommand{\eps}{\epsilon}

\newcommand{\bt}{\mathbf{t}}
\newcommand{\bs}{\mathbf{s}}
\newcommand{\bJ}{\mathbf{J}}

\newcommand{\EEnd}{\mathrm{EEnd}}

\begin{document}

\title[Gradient vector SRB Entropy]{Lipschitz Continuity and Formulas  of the Gradient Vector of the SRB Entropy Functional }

\author{Jianyu Chen${}^1$    and Miaohua Jiang${}^2$  }
\address{
${}^1$ School of Mathematical Sciences,
Center for Dynamical Systems and Differential Equations, Soochow University, Suzhou, Jiangsu, P.R.China.
}
\address{
${}^2$ Department of Mathematics, Wake Forest University, Winston Salem, NC 27109, USA}

\email{ ${}^1$jychen@suda.edu.cn;
\medskip
${}^2$jiangm@wfu.edu}

%\thanks{    }
%\subjclass[2020]{Primary: 37E05; Secondary: 37A35, 37A60, 80-10, 82B05}

\begin{abstract}  The Sobolev embedding theorem implies that the SRB entropy functional is also differentiable in the family of Anosov diffeomorphisms equipped with a suitable Hilbert manifold structure. The same holds true for the SRB entropy functional over the family of smooth expanding maps on a closed Riemannian manifold.  This implication leads to the local existence of the gradient flow of the SRB entropy and an explicit formula of the gradient vector of the entropy functional via the linear response of the SRB measure.
\end{abstract}

\maketitle

\large
\baselineskip 16 pt

\centerline{\bf\Large }

\section{Introduction }

The goal of the study is to investigate properties of  the Kolmogorov-Sinai (a.k.a. measure-theoretic or metric) entropy with respect to the Sinai-Ruelle-Bowen measure (SRB entropy for short)  as a functional on families of uniformly  hyperbolic systems, including expanding maps \cite{Young}. The motivation is to understand  whether  the SRB entropy of chaotic dynamical systems can be regarded as the corresponding quantity of the Boltzmann entropy of classical thermodynamics under the Gallavotti-Cohen Chaotic Hypothesis \cite{J21}.

In this article, we address two basic questions concerning the SRB entropy functional: (1) Does the entropy functional define a gradient flow in the space of uniformly hyperbolic systems equipped with a suitable Hilbert manifold structure?  (2)  If it does, can we derive the explicit formula of the gradient vector? The answers to both questions are affirmative.

  We limit our discussion to  two families of uniformly hyperbolic systems,
     transitive Anosov maps  and locally expanding maps,  defined on a closed Riemannian manifold such as an $n$-dimensional torus, extending some results in \cite{J24} for low dimensional systems to systems of any finite dimension.

\section{Sobolev Manifold Structure in the Space of Diffeomorphisms or Endomorphisms of a Closed Riemannian Manifold   }

It is now well-known that the SRB entropy of a uniformly hyperbolic system depends on \emph{the system} differentiably when the family of the systems is equipped with a suitable $C^r$ topology, $r \ge 3$.
However, it is preferable to have a Hilbert manifold structure for the family  so that the gradient vector field of the
functional is well-defined.

Since the Riemannian manifold $M$ is assumed to be compact, using the exponential map on $M$, it can be covered by a finite family of coordinate charts $\{(U_i, \phi_i)\}_{i=1}^l$, where $\phi_i: U_i \subset  M \to \R^n$ is a diffeomorphism for each $i$.  Indeed,  we may consider each $U_i$ an open ball of $\R^n$. Recall that a Hilbert (or more generally, Banach)  manifold is a topological space $\mathcal M$   such that
every point  in $ \mathcal M$ has a neighborhood homeomorphic to an open subset of a Hilbert (Banach) space.
The collection of such homeomorphisms (charts) forms an atlas, and the transition maps between overlapping charts are smooth (i.e., $C^r$ Fr\'echet differentiable for some $r \geq 1$).
For any diffeomorphism of $M$, its small $C^r$ open neighborhood is identified with an open neighborhood of the zero section of the tangent bundle via the exponential map of $M$, i.e, a $C^r$ vector field over $M$. Thus,
 to define a Hilbert manifold structure for the family of diffeomorphisms of $M$,  we just need to define a Sobolev norm for the space of differentiable vector fields over an open set $U \subset \R^n$. This Sobolev space will help to define the tangent space at any given diffeomorphism or endomorphism of $M$.

Let $V^{C^\infty}(U, \R^n)$ denote the family of $C^\infty$  vector fields defined on $U \subset \R^n$. Each vector field $F \in V^{C^\infty}(U, \R^n)$ can be described by its coordinate functions: $F(x)= (F_1(x), \cdots, F_n(x))$, where $F_i: U \to \R$ are $C^\infty$ functions. Let $\partial^\alpha F_i $ denote any partial derivative of $F_i$ of the order $\alpha=(\alpha_1, \cdots, \alpha_n)$.  The Sobolev norm $\|\cdot \|_{H^{k,p}} $ of $F$   is defined by
\begin{equation}\label{eq:norm}
  \| F\|^p_{H^{k,p}} = \sum_{i=1}^n  \sum_{0\le |\alpha| \le  k}  \int_{U}  |\partial^\alpha  F_i |^p dx ,\end{equation}
where $p>1, |\alpha|= \alpha_1 + \cdots +\alpha_n $ and the integral is just the Lebesgue integral.

Throughout the paper, $p$ is set to be $2$.  Then,  $\| F\|_{H^{k,2}} $ is a Hilbert  norm  on $V^{C^\infty}(U, \R^n)$. We omit the second superscript $2$ for simplicity.

The Sobolev space $H^k(U, \R^n)$ denotes the completion of the space $V^{C^\infty}(U, \R^n)$
%$\{   F:  C^\infty\text{vector\  field\ over }  U  \}$
under the norm $ \| F\|_{H^k} $. It is a Hilbert space \cite{He}.

Once this Hilbert space is defined over each coordinate chart, we put them together using a partition of unity $\psi_j(x), j =1,2, \cdots, l$ subordinated to the atlas $\{U_j\}_{j=1}^l$ to define a Sobolev norm for any diffeomorphism $f$ of $M$:
\begin{equation}\label{eq:norm2}
  \| f\|^2_{H^{k}} = \sum_{j=1}^l \sum_{i=1}^n  \sum_{0\le |\alpha| \le  k}  \int_{U_j} \psi_j(x)    |\partial^\alpha  f_i |^2 dx,\end{equation}
  where $f_i$ denotes the $i$th component of $f$ over the chart $U_j$ when it is identified with a vector field over $U_j$.

 Since we only need the distance between two nearby diffeomorphisms $f, g$, they can be identified as two vector fields over each coordinate chart. Thus, the distance between $f $ and $g$ is simply $\|f -g \|_{H^k}.$ The inner product between $f$ and $g$ is
 defined accordingly consistent with the definition of the Hilbert norm (\ref{eq:norm2}).

Other spaces we will use in the paper can be defined similarly and can be found in books such as  \cite{Mane, KH}.
We list them below for convenience.

$ \circ\  C^r(M), r \ge 0 $:  Banach space of $C^r$ functions over $M$. $r\ge 0$ can be any real number.

$\circ\  H^k(M), k \ge 0 $:  Hilbert space that is the completing of $C^k$ functions over $M$ under the Sobolev norm of and integer order $k\ge 0$.

We note that $C^k(M) \subset H^k(M)$ for any integer $k \ge 0$. There exists a constant $C$ such that $\|\varphi\|_{H^k} \le C  \|\varphi\|_{C^k} $
for any $\varphi \in C^k(M) $.

$\circ\  \diff^{C^r}(M), r \ge 0 $:  Banach manifold of $C^r$ diffeomorphism over $M$. $r\ge 0$ can be any real number.

$\circ\  \diff^{H^k}(M), k \ge 0 $:  Hilbert manifold that is the completion of $C^k$ diffeomorphisms over $M$ under the Sobolev norm of an integer order $k\ge 0$.  Note also $ \diff^{C^k}(M) \subset \diff^{H^k}(M).$

Similarly, we define $ \Endo^{C^r}(M), r \ge 0 $, $ \Endo^{H^k}(M), k \ge 0 $ for endomorphisms on $M$.

\begin{comment}
{\color{red}
Maybe recall the definition of Banach manifold and Hilbert manifold somewhere ??
Especially the definition of inner product on $\diff^{H^k}(M)$ and $\diff^{H^k}(M)$ ??
} - Revised: see Page 2.
\end{comment}

For $k\ge 1$, the family of transitive Anosov maps in $\text{Diff}^{H^k}(M)$ will be denoted by $A^{H^k}(M)$ and
the family of  expanding maps in $\text{Endo}^{H^k}(M)$ will be denoted by $E^{H^k}(M)$. Both families are Hilbert manifolds  with an inherited  Hilbert manifold structure.  Spaces $A^{C^r}(M)$ and $E^{C^r}(M)$ are also defined similarly.

Based on a Sobolev space embedding theorem that we will state in next section,  for any two integers $k,r$ satisfying  $k-r > \frac{n}{2}$ where $n$ is the dimension of the manifold $M$,  any vector field in
$H^k(U, \R^n)$ is also $C^r$, i.e., $ H^k(U, \R^n) \subset  C^r(U, \R^n)$ . Thus,
when $k > 2 +  \frac{n}{2}$, the SRB measure $\rho_f$ exists uniquely for any transitive Anosov map $f  \in
 A^{H^k}(M)$ (or any expanding map in $E^{H^k}(M)$).  The metric entropy of $f$ with respect to $\rho_f$ defines a functional
 from the Hilbert manifold $A^{H^k}(M)$ (or, $E^{H^k}(M)$) to $\R.$ We denote this entropy functional by $\CH(f)$.

\section{ Differentiability of the SRB Entropy in the Space  $A^{H^k}(M)$ or $E^{H^k}(M)$ }

The following embedding theorem \cite{He}  guarantees that the SRB entropy $\CH(f)$  is a differentiable functional
on Hilbert manifolds    $A^{H^k}(M)$ and $E^{H^k}(M)$ when $k > 3 + \frac{n}{2}$.

\begin{comment}
{\color{red}
The notation for Sobolev space is changed to $H^{k, p}(M)$ instead of $H^p_k(M)$ ??
} - OK, will check consistency.
\end{comment}

\begin{ebtheorem}\label{embedding} \cite{He} \text{\rm (Sobolev Space Embedding Theorem) }
  For any real number   $p \ge 1$ and two integers $  k,m$, if $k -m > \frac{n}{p}$, then $H^{k, p} (M) \subset C^m(M)$,
  where $M$ is an $n$-dimensional closed Riemannian manifold, $C^m(M)$ is the family of functions on $M$ with continuous derivatives up to order $m$.
\end{ebtheorem}

This embedding theorem implies that we have
$\text{Diff}^{H^k}(M) \subset  \text{Diff}^{C^m}(M)$
and
 $A^{H^k}(M) \subset    A^{C^m}(M)$
when $k -m >   \frac{n}{2}$.
The same holds for the space of endomorphisms.

\textbf{Remarks: }   (1)  The embedding of $H^{k, p} (M)$ in $C^m(M)$ is understood in the following sense: for each map $f$ in $H^{k, p} (M)$, we can change its values over a Lebesgue measure zero set such that $ f \in C^m(M)$ since maps in $H^{k, p} (M)$ are defined by equivalent classes. Two maps in  $H^{k, p} (M)$  are considered the same if their values differ only over a measure zero set.

(2)   We also need the local version of this embedding theorem: if $U\subset M$ is an open subset satisfying the cone condition (for example, an open ball), then the same embedding theorem holds:
$H^{k, p} (U) \subset C^m(U)$. Indeed, the proof of Theorem \cite{He} is carried out over the finite coordinate charts of $M$. See Pages 34-35 of Hebey's book \cite{He} and Pages 79-85 Adams and Fournier's book \cite{AF}.

We provide a proof for the differentiability of the SRB entropy functional $\CH(f)$  in the Sobolev norm.   Let $\rho_f$ be the unique SRB measure defined for
$f \in A^{H^k}(M)$ or $ E^{H^k}(M)$ and ${\mathcal H}(f)$ denote the Kolmogorov-Sinai  entropy  of $f$ with respect to $\rho_f$.

\begin{theorem}\label{differentiable} The SRB entropy functional
${\mathcal H}(f)$ is Fr\'echet differentiable on $A^{H^k}(M)$ (or  $ E^{H^k}(M)$), where  $k > 3 + \frac{n}{2}$ and $n$ is the dimension of the Riemannian manifold $M$.
\end{theorem}

\begin{proof}
  We know that when $k > 3 + \frac{n}{2}$,  the functional $\CH(f)$ is a Fr\'echet differentiable functional on
  $A^{C^3}(M)$ \cite{Ru97, Ru03}. Let's denote the derivative operator at the point $f$ by $D\CH_f$. We have
  \[  \lim_{\epsilon \to 0} \frac{1}{\epsilon}  |  \CH(f +\epsilon g) -  \CH(f ) - D\CH_f g | = 0 ,\] for all $g $ in the tangent space $T_f A^{C^3}(M)$.   Since  $A^{H^k}(M) \subset  A^{C^3}(M)$, the tangent space of $A^{H^k}(M)$ is a subspace of the tangent space of  $A^{C^3}(M)$ at $f$. Thus, the derivative operator $D\CH_f$ is a linear operator  defined on the tangent space of $A^{H^k}(M)$.

  Since $D \CH_f$ is a bounded operator on the tangent space $T_f A^{C^3}(M)$ and there exists a constant $C_2$ such that $\|g\|_{C^3} \le C_2 \|g\|_{H^k}$ for all $g \in T_f A^{C^3}(M)$ by the embedding theorem, we have
  \[   | D\CH_f g   |  \le  C_1 \|g\|_{C^3} \le  C_1 C_2 \|g\|_{H^k},\]  for some constants $C_1$ and $C_2$. Thus,
  $D\CH_f$ is a bounded linear operator on the tangent space of $A^{H^k}(M)$ at the point $f$,  and the SRB entropy functional  $\CH(f)$ is Fr\'echet differentiable on $A^{H^k}(M)$.
\end{proof}

An immediate consequence of this theorem is that at each point $f \in A^{H^k}(M)$, the gradient vector of $\CH $ exists uniquely since $A^{H^k}(M)$ is a Hilbert manifold.   Recall that the gradient vector at $f$, denoted by $\nabla \CH_f  $ is defined to be the unique vector in the tangent space $T_f A^{H^k}(M)$ such that   $<\nabla \CH_f, g> =  D\CH_f g $ for every $g \in T_f A^{H^k}(M)$.
Thus, we may also denote this gradient vector by $\nabla \CH_f =D\CH_f$, when $f \in  A^{H^k}(M)$.

We now consider whether this gradient vector field $ \nabla \CH  $ is integrable on the Hilbert manifold $A^{H^k}(M)$ .

Recall that a vector field $F$ defined in an open set $U$ of a Banach space $B$ is called locally integrable if for every point $x \in U$, there exist an $\epsilon>0$ and a map $\Phi(t,x)$ from $\   (-\epsilon, \epsilon) \times U \to U$ differentiable in $t$ such that $\Phi_t(x) :=\Phi(t, x),     U \to B  $   satisfies $\Phi_0(x)=x, \Phi_t( \Phi_s(x))= \Phi_{t+s}(x),  s,t, s+t \in   (-\epsilon, \epsilon),$ and $\frac{d}{dt}\big|_{t =0}  \Phi_t(x) = F(x).$

\begin{comment}
{\color{red}
Should be $\dots \to B$?  - corrected.
}\end{comment}

If  a vector field $F(x)$ is Lipschitz continuous in $x$, then, it is locally integrable \cite{Gi}, i.e., the initial value problem
$$ \dot x(t) = F(x),  x(0)=x ,$$ has a unique solution $\Phi_s(x)$ for $s\in   (-\epsilon, \epsilon):$  $\frac{d}{dt}\big|_{t =s}  \Phi_t(x) = F(    \Phi_s(x)).$

  We recall that the derivative operator $D\CH_f$ is  $C^{m-2}$ in $f$  when $m \ge 3$ \cite{Ru97}.  However, this does not imply $\nabla \CH_f$ is differentiable in $f$ with respect to the Sobolev norm since the derivative operator $D\CH_f$ is a linear operator on  the tangent space $T_fA^{C^m}(M)$ and thus,  belongs to the dual space of $T_fA^{C^m}(M)$,  which is a proper subspace of the dual space of $T_f A^{H^k}(M)$.

  On the other hand, we only need the Lipschitz continuity of the gradient vector field  $\nabla \CH_f$ in $f \in A^{H^k}(M)  \subset A^{C^m}(M)$ in the Sobolev norm. It follows from a direct estimate: For any given maps $f_1, f_2 $ close in the Hilbert manifold $ A^{H^k}(M) $, we may assume that they are in the same coordinate chart. Indeed, due to the structural stability of Anosov maps, we can assume $f_2$ is in an $H^k$ neighborhood of $f_1$.  Thus, $D\CH_{f_1} $ and $D\CH_{f_2} $ are acting on the same tangent space.

We have for any $g \in T_{f_1}  A^{H^k}(M) $,
  \begin{equation*}
  \begin{split}
  \left| <(\nabla \CH_{f_1} - \nabla \CH_{f_2}),  g> \right|
  &=
  \left| (D\CH_{f_1}  - D\CH_{f_2}) g  \right|
  \le C_1 \| f_1 - f_2 \|_{C^m}
  {  \|g\|_{C^m}}\\
  & \le C_1 C_2^2 \| f_1 - f_2 \|_{H^k} {    \|g\|_{H^k}}.
 \end{split}
  \end{equation*}

We have the following conclusion for both  Hilbert manifolds $A^{H^k}(M)$ and  $E^{H^k}(M)$.

\begin{theorem}\label{localexist}
  The gradient vector field $ \nabla \CH $  of the  SRB entropy $\CH(f)$ over the Hilbert manifold $A^{H^k}(M)$ (or $E^{H^k}(M)$)   is locally integrable: for each $f \in  A^{H^k}(M)$ (or $E^{H^k}(M)$), there exist  an $\epsilon_0$ and a unique map  $\Phi_t( f)=\Phi(t, f),  (t, f)\in  (-\epsilon_0, \epsilon_0) \times U_{\epsilon_0}(f_0)$ such that the map is differentiable in $t$ and satisfies $\Phi(0, f) = f,  \Phi(t+s, f ) = \Phi(t, \Phi(s, f)) ,$ when $ s,t, s+t \in (-\epsilon_0, \epsilon_0)$ and $f \in U_{\epsilon_0}(f_0)$, where $U_{\epsilon_0}(f_0)$ is the $\epsilon_0$-neighborhood of $f_0$ in $A^{H^k}(M)$ (or $E^{H^k}(M)$).
\end{theorem}

{\bf Remark }
While the proof of the differentiability of $\CH(f)$ over $A^{C^m}(M),  m\ge 3$ via thermodynamic formalism can be found in \cite{Ru97,Ru03}.    The proof of the differentiability of $\CH(f)$ over $E^{C^m}(M), m\ge 3$ is quite different: it is a consequence of the perturbation theory for transfer operators developed by Keller and Liverani \cite{KL}. We will provide more details in Appendix \ref{app: linear response}.

\section{Questions Arising from the Gallavotti-Cohen Chaotic Hypothesis }

By the uniqueness of the  local flow, given an initial map $f_0 \in A^{H^k}(M)$ or $E^{H^k}(M)$ , the orbit of the flow $\Phi_t(f_0)$ starting from $f_0$ can be extended in both directions of time $t$.

Natural questions arise concerning the global behavior of the flow governed by the gradient vector field $\nabla \CH$, in particular, in view of Gallavotti-Cohen's Chaotic Hypothesis \cite{GC, G96, G06} and  Maximum Entropy Production Principle \cite{JP, Maas},  the gradient flow $\Phi_t(f)$ may be regarded as a mathematical model for the process of a thermodynamic system evolving to its equilibrium.

The first question concerns the existence of a global equilibrium. The basic postulate of thermodynamics states that for a closed thermodynamic system, an equilibrium exists: any system not at equilibrium will evolve to a unique equilibrium.    The spirit of   Maximum Entropy Production Principle \cite{JP}  suggests that a system not at equilibrium will evolve in the direction that the entropy production is maximum, i.e., in the direction of the gradient vector of the entropy functional.   We can now formulate the corresponding statements as conjectures in the context of families of Anosov systems $A^{H^k}(M)$  and  expanding maps   $E^{H^k}(M)$.
Any map $f \in A^{H^k}(M)$ or  $E^{H^k}(M)$ is called at equilibrium if $\nabla \CH_f =0$.

{Notice that in any path-connect component of  $A^{H^k}(M)$ or  $E^{H^k}(M)$, all maps are topologically conjugate to each other due  to the structural stability. Thus, $\CH(f)$ has an upper bound equal to the topological entropy  of any map $f$ in the path-connected component, denoted as
$h_{\text top}(f)$.
}

\begin{comment}
{\color{red} The above statement need the following restriction:
in a topological conjugacy class $A^{H^k}_g(M)$
of $A^{H^k}(M)$,
i.e., all $f\in A^{H^k}(M)$ which is topologically conjugate to $g$,
the global maximum of $\CH(f)$ in $A^{H^k}_g(M)$
is equal to $h_{\text top}(g)$.
In general, there is no global maximum of $\CH(f)$ in $A^{H^k}(M)$,
since if $g\in A^{H^k}(M)$, then $g^n\in A^{H^k}(M)$ for any $n\ge 1$,
and $\CH(g^n)=n\CH(g)$ goes to $+\infty$.
}
\end{comment}

\emph{
Conjecture 1:}  (Global integrability of the gradient flow)   The gradient flow $\Phi_t( f)$ exists for all $t \in (-\infty, \infty)$ for every $ f \in  A^{H^k}(M)$ or  $E^{H^k}(M)$.

\emph{
Conjecture 2:}  (Existence of an equilibrium)  Given any $f \in A^{H^k}(M)$ or  $E^{H^k}(M)$ with $\CH(f) < h_{\text top}$, $\lim_{ t \to  \infty}  \Phi_t( f)$ exists.

\begin{comment}
{\color{red} Some thoughts: I think for the global equilibrium,
we need to first figure out if the maximum of SRB entropy could appear in
the border of $A^{H^k}(M)$. If it does, a solution of gradient flow may evolve to
a non-Anosov equilibrium.
My worry is that  under the evolution of SRB entropy gradient flow,
there might be a solution as a family of volume preserving Anosov diffeomorphisms on $\TT^3$,
which has two-dimensional unstable and one-dimensional stable,
such that  the limit becomes not Anosov, but partially hyperbolic (i.e., a central direction appears),
since the SRB gradient does not guarantee uniform hyperbolicity.
Anyhow, in  low dimensional case, say $\TT^2$, the conjecture of global equilibrium seems
doable, by delicate estimation of the sizes of stable/unstable cone families and their angle
under the gradient flow to obtain uniform hyperbolicity.

MJ- Agreed: one hopes that when a map evolves along a gradient trajectory, the optimal increase of the SRB entropy will force the hyperbolicity to become  more uniform in the phase space. In the toral Anosov map case, if a path-connected component contains an automorphism, the global equilibrium will be an Anosov map smoothly conjugate to the automorphism.
}

\end{comment}

If both conjectures hold true, we may say that  the basic postulate and the second law of thermodynamics are realized in these two mathematical models of thermodynamics. They give supporting evidences to Gallavotti-Cohen Chaotic Hypothesis and  Maximum Entropy Production Principle for closed nonequilibrium systems.

So far, progresses in proving these two conjectures are only made in special cases when the manifold $M$ is of low dimension.

\begin{comment}
Note: Conjecture 2 should be true for expanding maps in multi-dimensional case. See J. C's notes: \emph{SRB entropy (6-5-2025).pdf}.
\end{comment}

In \cite{J21, JL22},  partial results were proven for expanding maps on a unit circle and Markov transformations on a closed interval.

In \cite{J24}, the global existence ($t \in [0, \infty)$) of the gradient flow  is established for the family of measure-preserving expanding maps on the unit circle $\T^1$. The proof utilizes the derivative formula of the SRB entropy functional and the Riesz representation of the gradient vector.  While the proof of the differentiability of the SRB entropy requires a higher order differentiability of the map $f$,  this  (higher order differentiability) condition may be weakened if we can have an explicit formula of the SRB entropy. For example, in the case of measure-preserving circle map case,  the entropy functional $\CH(f)$ is given by  $\CH(f)= \int_{\T^1} \log f'(x) dx .$   This functional is differentiable in $f$ when $f$ is considered in the  family of $C^{1+\alpha}$ maps. With a lower order differentiability requirement, it is possible to obtain the Riesz representation of the gradient vector in a simple form, which leads to the proof of the global existence and to numerical approximations of the gradient flow.

In next two sections, we will give explicit formulas that characterize the gradient vector $\nabla \CH$ at any give point $f$.  The formula in the case of Anosov maps was essentially derived by Ruelle in \cite{Ru97,Ru03}. We restate the formula in the present context of the gradient vector of the SRB entropy.  The formula in the case of expanding maps on a closed Riemannian manifold of any finite dimension is new. We provide the derivation of the formula, including differentiability of the SRB entropy functional in full detail in the last section.

\section{Formulas characterizing the  Gradient Vector of the SRB Entropy}

Since the SRB entropy functional $\CH(f)$ is Fr\`echet differentiable on the Hilbert manifold $A^{H^k}\!(M)$ and
$E^{H^k}\!(M)$, its gradient vector at any given $f$, $\nabla H_f$ is a vector in the tangent space $T_f A^{H^k}(M)$, which is identified with the Hilbert space $V^{H^k}(M)$, the space of $H^k$ vector space over $M$. Thus, to characterize the gradient vector $\nabla H_f$ , we need to determine the value of
$\left< \nabla H_f, X\right> $ for any $H^k$ vector field $X\in V^{H^k}(M)$.

\subsection{ Gradient vector formula in the Anosov case}

Given that $f$ is $C^{1+\alpha}, \alpha >0$, transitive Anosov map on a closed manifold $M$,
the entropy formula with respect to its unique SRB measure $\rho_f$ is given by \cite{Mane}:
$$\CH(f)=  \int_{M}  \log J^uf  d \rho_f .$$

Ruelle first derived its derivative formula following his calculation of the linear response function of the SRB measure of $f$: the derivative formula for the SRB measure $\rho_f$ with respect to $f$. For details, see \cite{Ru97, Ru03}    and \cite{J12}.

The formula is derived in the more general context of hyperbolic attractors. We assume now that $f$ is $C^3$ Anosov on a closed manifold $M$: $ f \in A^{C^3}(M)$ and state the version of the theorem that we need.  Any $C^3$ small change of $f$ can be uniquely represented as a $C^3$ vector field $\delta f \in V^{C^3}(M)$ evaluated at $f(x)$. Let  $X= \delta f \circ f^{-1}$ be the pullback of $\delta f$. $X=X^u+X^s$ denotes the projection of the vector field $X$ onto the unstable and stable subbundles of the tangent bundle invariant under the derivative operator $Df$.

 \begin{theorem}  The SRB entropy functional  $$f \to   \CH(f)=  \int_{M}  \log J^uf  d \rho_f  $$
is $C^1$ in a $C^3$ neighborhood of $f$  in the $C^3$-norm.
 The derivative formula of the entropy functional $\CH(f)$ in the direction of a vector field $X= \delta f \circ f^{-1}$  is given by
$$ D\CH(f) X = \sum_{k \in \mathbb{Z} }  \int_M       \Div^u_{\rho}  X^u    \log J^u_f \circ f^k  \   d    \rho_f,$$
where    $  \Div^u_{\rho}   X^u  $     is the divergence of the vector field $X^u$ over the unstable manifold  of $f$ with respect to the volume form defined by the density function of the conditional measure of the SRB measure $\rho_f$ and $J^u_f$ is the Jacobian of $f$ along the unstable manifold.
\end{theorem}

We provide a sketch of the proof. For details, see Ruelle \cite{Ru97,Ru03} and Jiang \cite{J12}.

\begin{proof}
We start with the entropy formula   $\CH(f)=  \int_{M}  \log J^uf  d \rho_f  $.
Due to the non-differentiability of the unstable subspaces with respect to the base point, the functional $\log J^u f$ is usually not differentiable in $f$. However, when $f$ is restricted to a small $C^3$ neighborhood of a given Anosov map $f_0$, $f $ is topologically conjugate to $f_0$ through a H\"older continuous homeomorphism $h_f$ on $M$:
$$ f \circ h_f = h_f \circ f_0.$$ Make a change of variables on $M$ provided by the map $h_f$. We have
$$\CH(f)=  \int_{M}  \log J^uf  d \rho_f  =   \int_M  \log J^u f \circ h_f  d \rho^*_f,$$
where $\rho^*_f $ is the measure obtained through the change of variables:
$$ \int_M  \varphi  d \rho_f = \int_M  \varphi \circ h_f  d \rho^*_f , \ \text{
i.e.} \    \rho_f(A) = \rho^*_f( h_f^{-1} (A)),$$ for every Borel subset  $A$ of $M$.

The functional  $  \log J^u f \circ h_f$ is now differentiable in $f$ in a small $C^3$ neighborhood of $f_0$. The measure
$\rho^*_f $, invariant under $f_0$, is the unique equilibrium state for the potential function $-  \log J^u f \circ h_f$ for $f_0$ and it is also differentiable in $f$ since an equilibrium state depends on the potential function differentiably based on a general result of the thermodynamic formalism \cite{Ru97}.

Furthermore, given a direction $X=\delta f \circ f^{-1}$ a $C^3$ vector field over $M$.
$$ D\CH(f) X  = \left<D_X (\log J^uf\circ h_f),  \rho^*_f\right> + \left< \log J^uf \circ h_f, D_X \rho^*_f\right>,$$
where we have dropped the subscript in $f_0$ for simplicity.
The first terms gives us \cite{Ru03,J12}
$$  \left<D_X (\log J^uf\circ h_f),  \rho^*_f\right>  =   \int_M   \Div^u_\rho  {X^u}  d \rho_f .$$
The second term of the derivative is given by Ruelle's linear response formula \cite{Ru97,Ru03}.  Combining two parts, we have
$$ D\CH(f) X = \sum_{k \in \mathbb{Z} }  \int_M       \Div^u_\rho  X^u\    \log J^u_f \circ f^k  \   d    \rho_f,$$
i.e., the correlation function of two H\"older continuous functions $ \Div^u_\rho  X^u $ and $\log J^uf $.
\end{proof}

By the Sobolev embedding theorem, we now have the following corollary. We  assume $f \in A^{H^k}(M)$, $ k > 4 + \frac{n}{2}$. A small perturbation of $f$ in $A^{H^k}(M)$ is uniquely represented as a vector field $\delta f \in V^{H^k}(M)$ evaluated at $f(x)$. Let  $X= \delta f \circ f^{-1}$ denote the pullback of $\delta f$ and $X=X^u+X^s$ the projection of $X$ unto unstable and stable invariant subspaces of the tangent bundle of $Df$.

 \begin{corollary}  The SRB entropy functional  $$f \to   \CH(f)=  \int_{M}  \log J^uf  d \rho_f  $$
is $C^1$ in $A^{H^k}(M)$ when $k > 4 +\frac{n}{2}$.
 The gradient vector field $\nabla\CH_f$ of the entropy functional $\CH(f)$ is characterized by
$$
\left<\nabla\CH_f,  X \right> = D\CH_f X = \sum_{k \in \mathbb{Z} }  \int_M       \Div^u_\rho  X^u    \log J^u_f \circ f^k  \   d    \rho_f,$$
where    $X \in H^k$,  $  \Div^u_\rho  X^u  $  is the divergence of $X^u$ along the unstable manifold of $f$ with respect to the volume form defined by the density function of the conditional measure of the SRB measure $\rho_f$ and $J^u_f$ is the Jacobian of $f$ along the unstable manifold.
\end{corollary}

 \textbf{ Remark}    The condition $k > 4 +\frac{n}{2}$ is not  necessary. Our main interest in this paper is to show the existence of the gradient flow and obtain an explicit formula that characterizes the gradient vector, not to determine the precise order of the Sobolev norm.

\subsection{The Expanding map case}

In the case of expanding maps on a closed Riemannian manifold, the statements and their proofs are different because the maps are not invertible.
Recall that $E^{C^r}(M)$  denotes the family  of $C^r, r \ge 1,$ locally expanding maps on a closed Riemannian manifold $M$ of dimension $n$. Since we are considering a flow in this space, we assume that
$f \in E^{C^r}(M)$  if $f$ is $C^r$ and all eigenvalues of its derivative  operator $Df(x)$ are located outside of the unit circle in the complex plane.
Since the universal covering of $M$ is necessarily  topologically
equivalent to  $\R^n$,
for simplicity,  we { only consider}
  $M=\T^n=\R^n/\Z^n$, the dimension $n$ torus equipped with a common Euclidean metric. \begin{comment}
{\color{red} Given an integer $b\ge 2$,
we further denote $E^{C^r}_b(M)$ the subspace of $E^{C^r}(M)$ such that
every $f\in E^{C^r}_b(M)$ }
has the same degree $b$, i.e., the number of preimages of $x \in M$ under $f$.\end{comment}
Note that we can also define $E^{C^{r}}(M) $ to be the space of expanding $C^r$ endomorphisms of $M$ when $r \in [0,1)$: there exist constants $\mu > 1$ and $\epsilon > 0$
such that $d(f(x), f(y)) \ge \mu d(x,y)$ for all $x,y \in M$ with $d(x,y) < \epsilon$. For each real number $r\ge 0$,  the $C^r$ norm $\| f \|_{C^r}$ is defined as usual \cite{Mane}.  For any $r_2 > r_1 \ge 0$, $E^{C^{r_2}}(M)  \subset E^{C^{r_1}}(M).$

When $r >  1$, each map $f$ has a unique SRB measure with a H\"older continuous positive density function $\rho_f$, which is the equilibrium state for the potential function $- \log Jf $, the negative logarithm of the Jacobian of  $f$.  The metric entropy of $f$ with respect to $\rho_f$ is again denoted by $\CH(f)$. It is given by the formula \cite{Mane}
$$\CH(f) = \int_{\T^n} \log Jf  \rho_f dx,$$ where the integral is with respect to the Lebesgue measure.

The family  $E^{C^r}(M)$ is a Banach manifold whose tangent space at each point $f$ is just the Banach space of all $C^r$ vector fields on $M=\T^n$.
A small $C^r$-neighborhood of any given map $f$ is thus, identified with an open $C^r$ neighborhood of the zero section of the tangent bundle over $\T^n$,i.e., an open $C^r$ neighborhood of the zero vector field over $\T^n$.
The space of the perturbation is denoted by $\CB^{r}$. It can be identified with the space of vector fields over $\R^n$ that are periodic in every component.  Given a small perturbation $g \in \CB^r  $ , $f + g $ denotes the map $ x \to  {\rm Exp}_{f(x)} g( f(x))$, where ${\rm Exp}_{f(x)}$ is the exponential map of $\T^n$ at the point $f$.  For simplicity, we can also use the vector space structure of $\R^n$ to define $f+g$:   Let $\tilde f$ and $\tilde g$ be the lifts of $f$ and $g$ to the universal covering $\R^n$, $f + g$  is precisely the projection of  $\tilde f + \tilde g$ from the universal covering to $\T^n$. The perturbed map $f + t g$ is defined in the same way for any sufficiently small number $t$.

For any $f \in E^{C^r}(M), r \ge 1$, we denote by $\CL_f$ the transfer operator induced by $f$:
\[\CL_f \varphi (x) = \sum_{y: f(y)=x} \frac{ \varphi(y)}{Jf(y)}.\] The transfer operator is well-defined on the space  $C^\theta(M), \theta \ge 0$, of $C^\theta-$functions on $M$.
It is a bounded linear operator when $r \ge \theta +1$.  It is known that the SRB entropy $\CH(f)$ is a differentiable functional on $E^{C^r}(M)$  when $r \ge 3$. What we want to derive in this section is its derivative formula which characterizes the gradient vector of the entropy  functional $\CH(f)$ when  it is restricted to a Hilbert submanifold $E^{H^k}(M)  \subset E^{C^r}(M)$ where $k-r > \frac{n}{2}$. In the process of the derivation of the derivative formula, we also provide a detailed proof of its differentiability.

\begin{theorem}
\label{thm: entropy diff} On the Banach manifold  $E^{C^r}(M), r \ge 3$,
the Gateaux derivative of
the SRB entropy functional $\CH(f)$ at a given point (map)  $f\in E^{C^r}(M)$ in the direction of
$g=(g_1, g_2, \cdots, g_n) \in \CB^{r}$ is given by
\beqn
\displaystyle
 D\CH(f)  g =
 \int_{\T^n}   \left(  \trace((Df)^{-1}Dg)          \rho_f   - \sum_{n=0}^\infty \CL_f^n \Div[\CL_f (g \rho_f)]\cdot \log Jf(x) \right) dx,
\eeqn
where $Jf=|\det(Df)|$ is the Jacobian of $f$, $\rho_f$ is the SRB density function, and $\Div[\CL_f (g \rho_f)]$ is the divergence of the vector field
$(\CL_f (g_1 \rho_f), \CL_f (g_2 \rho_f), \cdots, \CL_f (g_n \rho_f))$.
Or, equivalently, by the duality definition of the transfer operator,
\beqn
 D\CH(f) g =    \int_{\T^n}
 \left( \trace((Df)^{-1}Dg) \rho_f  + \sum_{n=0}^\infty   \Div[\CL_f (g \rho_f)]\cdot \log Jf(x)\circ f^n  \right) dx.
 \eeqn
Moreover, the map $f\mapsto D\CH(f)$  from $E^{C^r}(M)$ to $\left(\CB^{r}\right)^*$, the dual space of $\CB^r$,
is Lipschitz.
\end{theorem}

\medskip
By the Sobolev embedding theorem, we have the following corollary when $M=\T^n$.
\begin{corollary}
The SRB entropy functional $\CH(f)$ is Fr\`echet differentiable
in the Sobolev space $E^{H^k}(\T^n)$  when $k> 3 +\frac{n}{2}$.
and the gradient vector $\nabla \CH_f$ of the entropy functional  $\CH(f)$ at a given point  $f\in E^{C^r}(M)$ is characterized by  the
formula  \beqn
\left<\nabla \CH_f, g\right>=
  \int_{\T^n}     \trace((Df)^{-1}Dg)          \rho_f   - \sum_{n=0}^\infty \CL_f^n \Div[\CL_f (g \rho_f)]\cdot \log Jf(x)  dx,
\eeqn
where
$g=(g_1, g_2, \cdots, g_n) \in V^{H^k}(\T^n)$,  $Jf=|\det(Df)|$ is the Jacobian of $f$, $\rho_f$ the SRB density function, and $\Div[\CL_f (g \rho_f)]$ is the divergence of the vector field\hfill\break
$(\CL_f (g_1 \rho_f), \CL_f (g_2 \rho_f), \cdots, \CL_f (g_n \rho_f))$.
Or, equivalently, in view of the integration by parts,
\beqn
\left<\nabla \CH_f, g \right>=    \int_{\T^n}   \trace((Df)^{-1}Dg) \rho_f  + \sum_{n=0}^\infty   \Div[\CL_f (g \rho_f)]\cdot \log Jf(x)\circ f^n  dx
 .\eeqn

 Furthermore, the gradient vector  $\nabla \CH_f$ is Lipschitz continuous in $f \in E^{H^k}(\T^n)$ and thus, the gradient flow of the entropy functional $\CH(f)$ exists locally in an open neighborhood of every $f \in E^{H^k}(\T^n)$.
\end{corollary}

{\bf Remarks}

   (1) When the vector field $g$ that defines the perturbation has the form $g = X\circ f$ for a vector field $X$ on $\T^n$,   the linear response formula for  the SRB measure of expanding maps  can also be found in \cite{Ba18, Belse}.
     The complete proof of the linear response formula that leads to the formula of the gradient vector of the SRB entropy is not available in literature.  Here we provide a complete derivation based on an abstract Banach space operator theory in the Appendix.

 (2) Again $k > 3 +\frac{n}{2}$ is not optimal. There is an incentive to keep $k$  as small as possible in $E^{H^k}(M)$: When $k$ is small, the Hilbert space $V^{H^k}(M)$ is larger, which might help to find the Riesz representation of the gradient vector.  It seems likely to have a smaller $k$ if one uses a direct approach from \cite{Ba18} to consider the spectral properties and stability of the transfer operator on Sobolev spaces.

The proof of Theorem \ref{thm: entropy diff} largely depends on the formula of the linear response function of the SRB measure for expanding maps. We present the formula below and provide a detailed proof in Appendix.

\medskip

Given  $f\in E^{C^r}(M), r\ge 3$ and  $g\in \CB^r$, there exists $\eps>0$ such that
$f_t:=f+t g\in E^{C^r}(M)$ for any $t\in (-\eps, \eps)$.
Let $\CL_t$ be the transfer operator
associated with $f_t$, and
let $\rho_t$ be the SRB density of $f_t=f+t g$.
For simplicity,  denote $\CL=\CL_0$ and $\rho=\rho_0$.

\begin{theorem}[Linear response formula]
\label{thm: linear response}
 The map $f_t \to   \rho_t  $ is Fr\`echet differentiable  from $E^{C^r}(M)$ to $C^0(M)$. There exists bounded linear operator on $\CB^r(M)$ to $C^0(M)$
\beq
\label{eq: lin res 0}
\xi=\xi_f(g):=\left. \p_t \rho_t \right|_{t=0}
=-\sum_{n=0}^\infty \CL^n \Div \left(\CL \left(g \rho\right)\right).
\eeq
such that for any $g \in \CB^r(M)$
\[  \lim_{t\to 0} \left\| \frac{1}{t} (\rho_t - \rho_0) - \xi_f(g)  \right\|_{C^0(\T^n)} = 0. \]
\end{theorem}

\subsubsection{Proof of Theorem~\ref{thm: entropy diff}}
\label{sec: proof}

By Theorem~\ref{thm: linear response},
the Gateaux derivative of the SRB entropy $\CH(f)$ at $f\in E^{C^r}(M)$ in the direction of
$g\in \CB^{r}$ is given by
\begin{align}
D\CH(f)g
&=\lim_{t\to 0} \frac{\CH(f+tg) - \CH(f)}{t}   \\
&=\left. \p_t  \right|_{t=0} \left( \int_{\T^n}  \log\left(J(f + tg) \right) \rho_t(x) dx \right)  \notag \\
&=\int  _{\T^n}  \trace((Df)^{-1}Dg)  \rho(x)  dx +\int_\TT \xi (x) \log(Jf(x)) dx  \notag \\
&=\int_{\T^n} \trace((Df)^{-1}Dg)  \rho(x) dx - \sum_{n=0}^\infty \int_\TT  \CL^n   \Div[\CL(g\rho)] \cdot \log(Jf(x)) dx \label{gradop1}\\
&= \int_{\T^n}  \trace((Df)^{-1}Dg)  \rho(x) dx + \sum_{n=0}^\infty \int_\TT  \Div[\CL(g\rho)]  \cdot \log Jf \circ f^n  dx \label{gradop2}
\end{align}

{\it     End of the proof of Theorem~\ref{thm: entropy diff}          }

{\bf Remark}
  The entropy functional $\CH(f)$ is in fact, Fr\'echet differentiable in $E^{C^r}(M)$. We can show the derivative operator $D\CH(f)$ is a bounded linear operator on $ \CB^{r}$ for every $ f \in E^{C^r}(M)$ and is Lipschitz continuous in $ E^{C^r}(M)$ and thus in $E^{H^{k}}(M)$.

{

\begin{theorem}
Let $r\ge 4$. The map
$\displaystyle
 f \mapsto D\CH(f) $ from $E^{C^r}(M) $ to $  (\CB^{r})^*
$
is locally Lipschitz continuous, i.e.,
for any $f\in E^{C^r}(M)$, there is a $C^r$-neighborhood $\CU_f$
and a constant $L_f>0$
such that
for any $f_1, f_2\in \CU_f$, and any $g\in \CB^r$,
\beqn
\left| D\CH(f_1)g - D\CH(f_2)g  \right|\le L_f \|f_1-f_2\|_{C^r} \|g\|_{C^r}.
\eeqn
\end{theorem}

\begin{proof}
We recall that in the proof of Theorem~\ref{thm: entropy diff},
for any $f\in E^{C^r}(M)$ and any $g\in \CB^{r}$,
\beqn
D\CH(f)g =
\int  _{\T^n}  \trace((Df)^{-1}Dg)(x)  \rho_f(x)  dx +\int_\TT [\xi_f(g)] (x) \log(Jf)(x) dx,
\eeqn
where $\rho_f$ is the SRB density of $f$,
and $\xi_f(g)$ is given by Formula~\eqref{eq: lin res 0}.
Since $r\ge 4$, it is easy to see that the mappings
$\displaystyle
f\mapsto \trace((Df)^{-1}Dg)
$
and
$\displaystyle
f\mapsto \log(Jf)
$
are $C^3$ smooth and certainly locally Lipschitz continuous.

It remains to show
$\displaystyle
f\mapsto \rho_f
$
and
$\displaystyle
f\mapsto  \xi_f(g)
$
are also locally Lipschitz continuous in the following sense:
there is a constant $C_f>0$ such that
for any $f_1, f_2\in \CU_f$ and any $g\in \CB^r$,
\beq
\label{Lip f12}
\begin{split}
& \left\|\rho_{f_1} - \rho_{f_2} \right\|_{C^0}\le C_f\|f_1-f_2\|_{C^r},
\\
&\left\|\xi_{f_1}(g) - \xi_{f_2}(g) \right\|_{C^0}\le C_f\|f_1-f_2\|_{C^r}\|g\|_{C^r}.
\end{split}
\eeq
To this end, we consider the two-parameter family
\beqn
\begin{split}
f_\bt
&=t_1 f_1+(1-t_1) f_2 + t_2 g  \\
&= f + (f_2-f) + t_1 (f_1-f_2) + t_2 g
\end{split}
\eeqn
where $\bt=(t_1, t_2)\in \bJ:=[0, 1]\times (-\eps, \eps)$.
Note that $f_\bt$ is of the form \eqref{def f_bt} in Appendix \ref{sec.abs3}
with $f_0=f_2-f$,
$g_1=f_1-f_2$ and $g_2=g$.
Let $\eps_f$ be given by Lemma~\ref{lem: sp gap uniform}.
As long as the value of $\eps$ and the size of $\CU_f$
are sufficiently small, we can ensure that $\|f_\bt - f\|_{C^r}<\eps_f$
for all $\bt\in \bJ$.
Let $\rho_\bt$ be the  SRB density of $f_\bt$, the by
Theorem~\ref{thm: der rho bt},
there is $C_f>0$ such that for all $\bt\in \bJ$,
\beqn
\|\rho_\bt\|_{C^0}\le C_f, \ \
\|\p_1 \rho_\bt\|_{C^0} \le
C_f\left\|f_1 - f_2\right\|_{C^r},  \ \
\|\p_1 \p_2 \rho_\bt\|_{C^0} \le C_f \left\|f_1 - f_2\right\|_{C^r}\left\|g\right\|_{C^r}.
\eeqn
 It follows that
\beqn
\|\rho_{f_1}-\rho_{f_2}\|_{C^0}
=\|\rho_{(1, 0)}-\rho_{(0, 0)}\|_{C^0}
= \left\| \left.\int_0^1 (\p_1 \rho_\bt) \right|_{\bt=(t_1, 0)} dt_1 \right\|_{C^0} \le C_f\|f_1-f_2\|_{C^r}.
\eeqn
and
\beqn
\begin{split}
\left\|\xi_{f_1}(g) - \xi_{f_2}(g) \right\|_{C^0}
=\left\|\left.\p_2\rho_{\bt}\right|_{\bt=(1, 0)}-\left.\p_2\rho_{\bt}\right|_{\bt=(0, 0)} \right\|_{C^0}
& =\left\| \left.\int_0^1 (\p_1 \p_2 \rho_\bt) \right|_{\bt=(t_1, 0)} dt_1 \right\|_{C^0} \\
&\le C_f \left\|f_1 - f_2\right\|_{C^r}\left\|g\right\|_{C^r}.
\end{split}
\eeqn
The proof of \eqref{Lip f12} and thus this theorem
are complete.
\end{proof}

}

The corollary is stated when $M=\T^n$.

\begin{corollary}
 The SRB entropy $\CH(f)$ is Fr\'chet differentiable on the Hilbert manifold
$E^{H^k}(\T^n) $ when $k > 3 +\frac{n}{2}$ and the gradient vector at $f\in E^{H^k}(\T^n) $, $\nabla H_f$ is locally Lipschitz in $f$ and  characterized by the following formula
\begin{equation} \label{gradvec1}
\left<\nabla H_f,  g\right> = \int_{\T^n} \trace((Df)^{-1}Dg)  \rho(x) dx - \sum_{n=0}^\infty \int_\TT  \CL^n   \Div[\CL(g\rho)] \cdot \log(Jf(x)) dx ,\\
\end{equation}
where $g $ is a vector field in  $V^{H^k}(\T^n)$,  $\rho$ is the density function of the SRB measure of $f$,
and $\CL$ is the transfer operator induced by $f$.
\end{corollary}

\begin{comment}
There seems to be a problem here: the gradient vector depends on the metric. How does the Sobolev metric play a role here?
What is the relation between the derivative formula of the entropy functional and the gradient vector?
\end{comment}
\newpage

\appendix

\section{Abstract Perturbation Theory and the derivation of the linear response formula}\label{sec.abs}

The derivation of the linear response formula for the SRB measure of an expanding map
is an application of an abstract perturbation theory of bounded linear operators on Banach spaces developed by Gouezel, Keller, and Liverani \cite{GL,KL}.
Since the precise version we need did not appear in literature, we provide the details in this appendix.

\subsection{Spectral  Decomposition for a bounded linear operator with a Spectral Gap}
Assume that  $\mathcal L $ is a bounded linear operator on a real Banach space $\mathcal B$ with a $(1, \eta)-$spectral gap, i.e., $1$ is a simple eigenvalue and ${\rm Spect}(\mathcal L) \setminus \{1 \} \subset \{ z \in \mathbb C:  |z| \le \eta < 1 \} $, where  ${\rm Spect}(\mathcal L)$ is the spectrum of $\mathcal L$.  Then, we have the following well known spectral decomposition theorem for the operator $\mathcal L$.

\begin{theorem}\label{decomposition}
There are two bounded linear operators $\mathcal P $ and $\mathcal N$ on $\mathcal B$ such that
$\mathcal  L =\mathcal  P +\mathcal  N$, where $\mathcal P$ is a one-dimensional projection, i.e.,
$\CP^2 =\CP$ and ${\rm dim}({\rm Im}(\mathcal P))=1$,  $\mathcal N$ has a spectral radius no greater than $\eta$, and $\CP\CN=\CN\CP=0$.
Furthermore,  the spectral projection $\mathcal P$ is given by a contour integral:
\begin{equation}\label{contour}
 \mathcal P = \frac{1}{2 \pi i} \oint_\gamma  R(z) dz,\end{equation}
where $\gamma$ is any simple closed curve in the region $V=\{ z \in C:  |z| > \eta\}$ with eigenvalue one in its interior
and $R(z)=(z - \mathcal L)^{-1}$ is the resolvent of $\CL$, a bounded linear operator on $\CB$ for all $z\not\in {\rm Spect}(\mathcal L)$.
\end{theorem}

The theorem follows from a more general spectral decomposition theorem that
 can be found in  \cite{Kato} in the context of linear operators on complex Banach spaces. For operators on real Banach spaces, we obtain the same theorem  by the standard technique of complexification.

\subsection{A one-parameter family of operators with a spectral gap on nested Banach spaces}

In order to apply the spectral decomposition theorem to transfer operators defined by expanding maps and obtain the linear response function, i.e., the derivative of the SRB measure with respect to a parameter, we need to extend   Theorem \ref{decomposition} to the case of linear operators defined on three nested Banach spaces
  $\mathcal B_2 \subset  \mathcal B_1\subset \mathcal B_0. $
  For any vector $\varphi_i \in \CB_m, m = 1,2$,
  $$\|  \varphi \|_{\CB_{i-1}}  \le \|  \varphi \|_{\CB_{i}} .$$ In other words, $\CB_m$ is a subspace of $\CB_{m-1}$ endowed with a stronger norm.

 Let $J =(-\epsilon, \epsilon)$ for some $\epsilon > 0$. Assume that for each $t\in J$,  $\CL_t$ is a bounded linear operator defined on all three Banach spaces $\CB_m, m=0,1,2$ and $\CL_t|_{\CB_m} = \CL_t|_{\CB_{m-1}}, m=1,2$. For convenience, we simply write $\CL_t$ instead of $\CL_t|_{\CB_m}$ when the underlying space is clear.

  We now make the following Assumption (H1) for this one-parameter family of operators.

  (H1): There is $\eta < 1$ such that for all $t \in J$, $\CL_t|_{\CB_m}, i=1,2$ has a  $(1,\eta)$-spectral gap.

  Note that Assumption (H1) is not made for $\CL_t|_{\CB_0}$. In the rest of this subsection, we only consider
  $\CL_t|_{\CB_m}, m=1,2$ when the spectrum is concerned.

  The resolvent of $\CL_t$,  given by
 \begin{equation} R_t(z) := (z - \CL_t)^{-1},  \end{equation}
 is a well-defined bounded linear operator on both $\CB_m, m=1, 2$ for any $z ,  |z| > \eta,  z \not=1.$

 By the spectral decomposition Theorem \ref{decomposition}, we have spectral projections $\CP_{t, m}$ and $\CN_{t,m}$ such that
 $\CL_t|_{\CB_m} =  \CP_{t, m} +  \CN_{t,m}.$ Or, simply,
 $$    \CL_t   =  \CP_{t } +  \CN_{t}       ,$$  since $\CP_{t,2}, \CN_{t,2}$ are just restrictions of $\CP_{t,1}, \CN_{t,1}$ to $\CB_2$.

 Thus, for this one-parameter family of operators, we have the corresponding spectral decomposition theorem. \begin{theorem}\label{decomposition2}
 Under Assumption (H1), for each $t \in J$, there are two bounded linear operators $\mathcal P_t $ and $\mathcal N_t$ on $\mathcal B_{1,2}$ such that
$\mathcal  L_t =\mathcal  P_t +\mathcal  N_t$, where $\mathcal P_t$ is a one-dimensional projection, i.e., $P_t^2 =P_t$ and ${\rm dim}({\rm Im}(\mathcal P_t))=1$,  $\mathcal N_t$ has a spectral radius no greater than $\eta$,  and $\CP_t \CN_t = \CN_t\CP_t =0$.
 \end{theorem}

\subsection{Continuity and differentiability assumptions of operators in one parameter}
\label{subsec: assump}

In the application of the spectral decomposition theorem to transfer operators defined by expanding maps, a one-parameter family of bounded linear operators $\CL_{t,i}$ does not depend on $t$ continuously in $\CB_m$.  But as an operator from   $\CB_m $ to $\CB_{m-1}, i =1,2,$  $\CL_t$ can be differentiable in $t$. To be precise, we make the following differentiability assumption for this one-parameter family of linear operators.

(H2):  For any $t \in J$, the limit
\begin{equation} \partial_t \CL_t := \lim_{s \to t } \frac{1}{s-t} ( \CL_s - \CL_t) \end{equation}
 exists as a bounded linear operator from $\CB_2$ to $\CB_1$.

We remark that the limit in (H2) is defined in the strong topology from $\CB_2$ to $\CB_1$, i.e., for any vector $\varphi \in \CB_2$, we have
\[  \lim_{s \to t } \left\|
 \left[\frac{1}{s-t} ( \CL_s - \CL_t)  - \partial_t\CL_t \right] \varphi
 \right\|_{\CB_1} =0.\]
Note that under the assumption (H2), the mapping $t \to \CL_t$ is continuous from $J$ to the Banach space of bounded linear operators $L(\CB_2, \CB_1)$ in the following sense: for any $\varphi \in \CB_2$,
$\lim_{s\to t} \left\|(\CL_s - \CL_t)\varphi \right\|_{\CB_1} =~0$.

Assume $\CL_t, t \in J$, satisfies (H1) and let $\gamma$ be the circle centered at $1$ in the complex plane with a radius $\frac{1-\eta}{2}$ and $R_t(z)$ be the resolvent. Assumption (H3) is about the continuity of the resolvent as an operator from $\CB_1$ to $\CB_0$.

(H3):  For any $t \in J$,
\[\lim_{s \to t} \sup_{z \in \gamma} \|R_s(z) - R_t(z)\|_{L(\CB_1, \CB_0)}=0.\]

It follows from Assumption (H3) that for any $\varphi \in \CB_1$, we have
\[  \lim_{s\to t} R_s(z) \varphi = R_t(z) \varphi,\] in which
$R_s(z) \varphi $ and $ R_t(z) \varphi$ are regarded as vectors in the bigger space $\CB_{0}$.

\begin{comment}
Note the change from the original version.
Question:  do we need it to be true for i=2?
\end{comment}

We stress that this limit may not hold in $\CB_1$. The verification of this condition for transfer operators is not at all trivial. Thanks to the perturbation theory developed by Gouezel-Keller-Liverani  (\cite{GL}, Section 8, \cite{KL}) (see also Baladi's book \cite{Ba18}, Chapter 2). Assumption (H3) can be obtained via a Lasota-Yorke estimate  of $\CL_t, t \in J$ uniform in $t$ on the pair $(\CB_1, \CB_0)$.
{
In fact,
the Gouezel-Keller-Liverani  perturbation theory for a family of abstract operators
is summarized in Appendix A.3 of Baladi's book \cite{Ba18},
which is based on the Lasota-Yorke type inequalities (A.2)-(A.7) therein.
Under such condition (for the particular case $N=1$),
for any $ \delta\in \left(0, \frac{1-\eta}{2}\right)$,
we denote
$\displaystyle
V_\delta:=\left\{ z\in \mathbb{C}: \ |z|\ge \eta+\delta,  \text{and} \ |z-1|\ge \delta\right\},
$
and note that the circle $\gamma\subset V_\delta$.
It is shown that (see Inequality (A.8) of Appendix A.3 of Baladi's book \cite{Ba18})
for any $t\in J$,
there exist $\eps=\eps_{t, \delta}>0$, $C=C_{t, \delta}>0$ and $\eta=\eta_{t, \delta}>0$ such that
for any $s\in J$ with $|s-t|<\eps$ and any $z\in V_\delta$,
\beqn
\|R_s(z) - R_t(z)\|_{L(\CB_1, \CB_0)}\le C |s-t|^\eta,
\eeqn
from which Assumption (H3)  follows.
}

\subsection{ Derivative formula of the spectral projection}

From now on, we assume that a one-parameter family of bounded linear operators $\CL_t, t \in J$,  on
a nested sequence of three Banach spaces $\CB_m, m =0,1,2$ satisfy all three assumptions (H1)-(H3). The goal of this subsection is to prove that the spectral projection
$\CP_t:  \CB_2 \to \CB_0 $ is differentiable with respect to $t$.
Moreover, we shall provide an explicit formula of the derivative of $\CP_t$.

\begin{comment}
Note:  $\CP_t:  \CB_2 \to \CB_2$ is changed to $\CP_t:  \CB_2 \to \CB_0 $
\end{comment}

\subsubsection{} {\bf Differentiability and the derivative formula of the Resolvent}

We begin with proving the differentiability of the resolvent operator $R_t(z)=(z-\CL_t)^{-1}$ when $z \in V\setminus \{1\}$,
where  $V=\{ z\in \mathbb C: |z|>\eta\}$.  The following lemma shows that $z \to  R_t(z)$ is a meromorphic function on $V$.

\begin{lemma}\label{merom} For any $t \in J$ and any $z \in V\setminus \{1\}$,
\[ R_t(z) = \frac{\CP_t}{z-1} + \CQ_t(z),\]
where \begin{equation}\label{Qt}
 \CQ_t(z) := R_t(z)(I - \CP_t) = (z- \CN_t)^{-1}(I-\CP_t)
 \end{equation}
 is a homomorphic function on $V$.
\end{lemma}

\begin{proof}
 By Theorem \ref{decomposition2}, we have
 \[   \CL_t\CP_t = \CP_t \CL_t = \CP_t  \ {\rm and \ }    \CL_t(I- \CP_t ) =    (  I-\CP_t) \CL_t= \CN_t.\]
 Hence, $  R_t(z) = R_t(z) \CP_t +  R_t(I-\CP_t)$, where
 \[ R_t(z)\CP(t) =(z- \CL_t)^{-1} \CP_t = (z-1)^{-1}\CP_t, \]
 since $(z- \CL_t)\CP_t = z \CP_t - \CP_t$ and
 \[ R_t(I-\CP_t) = (z-\CL_t)^{-1} (I- \CP_t) = (z- N_t)^{-1} (I -\CP_t) =:  \CQ_t(z),\]
 The second equality holds in the last equation because for any $\psi = (I- \CP_t)\varphi$, $\CP_t \psi =0$ and $\CL_t \psi = \CN_t \psi$.

It remains to show that $\CQ_t(z)$ is homomorphic in the region $V.$  Indeed, since $|z| >\eta$ is greater than the spectral radius of $\CN_t$, we have the expansion
\[ (z-\CN_t)^{-1} = \sum^\infty_{n=0} z^{-(n+1)} \CN_t^n,\]
which is homomorphic in $|z| > \eta$. Therefore, $\CQ_t(z)=(z-\CN_t)^{-1}(I-\CP_t)$ is also homomorphic in $|z| > \eta$.
\end{proof}

Recall that the contour $\gamma$ is a circle centered at 1 with a radius $\frac{1}{2} (1 - \eta)$. The following lemma shows the differentiability of the resolvent $R_t(z)$ with respect to $t$ for any $z \in \gamma$.

\begin{lemma}\label{dev.R}
For any $t \in J$ and any $z \in \gamma$, we have
\[   \partial_t R_t(z):= \lim_{s \to t} \frac{1}{s-t}(R_s(z)-R_t(z) ) =R_t(z)(\partial_t \CL_t)R_t(z) ,\]
which is a bounded linear operator from $\CB_2$ to $\CB_0$.
\end{lemma}

\begin{proof}
By Assumption (H3), for any $z \in \gamma$, the map $s \to R_s(z)$ is continuous from $s \in J$ to the space $L(\CB_1, \CB_0)$. Due to the following resolvent identity
\[  R_s(z)-R_t(z) = R_s(z)(\CL_s -\CL_t)R_t(z),\]
the lemma follows immediately from (H2): Given any $\varphi \in \CB_2$, $R_t(z) \varphi \in \CB_2$.
\begin{align*} &[ \frac{1}{s-t}( R_s(z)-R_t(z) ) - R_t(z)(\partial_t \CL_t)R_t(z)]\varphi\\
= & [ \frac{1}{s-t}R_s(z)(\CL_s-\CL_t)R_t(z)  - R_t(z)(\partial_t \CL_t)R_t(z)]\varphi \\
 =& [ \frac{1}{s-t} R_s(z)(\CL_s-\CL_t)R_t(z)  -   R_s(z)(\partial_t \CL_t)R_t(z)  +     R_s(z)(\partial_t \CL_t)R_t(z)      -            R_t(z)(\partial_t \CL_t)R_t(z)             ]\varphi \\
 = & R_s(z)[ \frac{1}{s-t}(\CL_s-\CL_t) -  (\partial_t \CL_t)]R_t(z)  \varphi  +  (R_s(z)- R_t(z))(\partial_t \CL_t)R_t(z)                  \varphi .\end{align*}
Taking the limit $s \to t$, we see \[ \lim_{s\to t} \|[ \frac{1}{s-t}( R_s(z)-R_t(z) ) - R_t(z)(\partial_t \CL_t)R_t(z)]\varphi\|_{\CB_0} =0 \] for any $\varphi \in \CB_2$.
\end{proof}

\subsubsection{}{\bf Differentiability and derivative formula of the spectral projection}

Recall that  $\CQ_t(z)$ is given by (\ref{Qt}).  Let
\begin{equation}\label{Q1}
\CQ_t:= \CQ_t(1)=(I-\CN_t)^{-1}(I-\CP_t)=\sum^\infty_{n=0} \CN^n_t (I-\CP_t) = \sum^\infty_{n=0} \CL^n_t (I-\CP_t).\end{equation}
The following lemma shows the differentiability of the spectral projection  $\CP_t$ as a linear operator from $\CB_2$ to $\CB_0$ with respect to $t$.

\begin{theorem}\label{dev.P} The map $ t \to \CP_t \in L(\CB_2, \CB_0)$ is differentiable in $t \in J$. Moreover, for any $t \in J$, we have
\begin{equation}\label{eq.dev.P}
\partial_t\CP_t = \CP_t(\partial_t \CL_t)\CQ_t  + \CQ_t(\partial_t \CL_t) \CP_t,
\end{equation} which is  bounded linear operator from $\CB_2$ to $\CB_0$.
\end{theorem}

\begin{proof}
By (\ref{contour})  and Lemma \ref{dev.R}, we have
\[\partial_t\CP_t = \frac{1}{2 \pi i} \oint_\gamma \partial_t R_t(z) dz
=\frac{1}{2 \pi i} \oint_\gamma R_t(z)(\partial_t \CL_t)R_t(z) dz, \]
which is a bounded linear operator from $\CB_2$ to $\CB_0$. By Lemma \ref{merom},
\[ R_t(z)(\partial_t \CL_t)R_t(z) = (\frac{\CP_t}{z-1} + \CQ_t(z))(\partial_t \CL_t) (\frac{\CP_t}{z-1} + \CQ_t(z)) \]
is a meromorphic function in $V=\{z: |z| >r\}$  with a pole at $1$ of order 2 and its order 1 term is
\[\frac{1}{z-1}[ \CP_t(\partial_t \CL_t)\CQ_t(1) + \CQ_t(1)   (\partial_t \CL_t) \CP_t ].   \]
By the residue theorem, we have
\[\partial_t\CP_t = \frac{1}{2 \pi i} \oint_\gamma R_t(z)(\partial_t \CL_t)R_t(z)  dz  =  \CP_t(\partial_t \CL_t)\CQ_t + \CQ_t    (\partial_t \CL_t) \CP_t.\]
\end{proof}

In the application to the calculation of the linear response formula of the SRB measure for expanding maps, we will need to consider the composition $(\partial_t \CP_t) \CP_t$. Next lemma gives a simple formula for this product.

\begin{lemma}  For any $t \in J$, we have
\begin{equation}\label{PPP} (\partial_t \CP_t) \CP_t = \sum_{n=0}^\infty \CL^n_t(\partial_t\CL_t)\CP_t.  \end{equation}
\end{lemma}

\begin{proof}
By Theorem \ref{dev.P} and the fact that  $\CP_t \CQ_t=\CQ_t\CP_t =0$, we have
\[  (\partial_t \CP_t) \CP_t =\CQ_t (\partial_t\CL_t)\CP_t
=   \sum^\infty_{n=0} \CL_t^n(I-\CP_t)(\partial_t\CL_t)\CP_t. \]
We now show $\CP_t( \partial_t\CL_t)\CP_t =0$. Taking derivative on both sides of $\CL_t \CP_t = \CP_t$, we have
\[(\partial_t\CL_t)\CP_t + \CL_t(\partial_t \CP_t)= \partial_t\CP_t.                         \]
Multiplying $\CP_t$ on both sides of the equation and using $\CP_t\CL_t=\CP_t$, we obtain $\CP_t( \partial_t\CL_t)\CP_t =0$. The formula (\ref{PPP}) follows.
\end{proof}

\subsection{Linear Response Formula of the SRB measure for expanding maps}
\label{app: linear response}

We consider the following situation:

$f$ is any given $C^r, r\ge 3 $ expanding map on  a closed Riemannian manifold $M$. For simplicity, we assume $M$ is just a dimension $n$ torus $\mathbb T^n$.

\begin{comment}
We need $f$ to be $C^3$ since we want the transfer operator to act on $\CB_2= C^2(M)$. This may not be optimal.
\end{comment}

The space of all such maps forms a Banach manifold whose tangent space $TM$  is a Banach space consisting of  all $C^3$ vector fields on $M$. By lifting $f$ to the universal covering of $M=\mathbb T^n$, the space of all $C^3$-vector fields on $M$ can be identified with the space of maps from $\R^n$ to itself that are $1$-periodic in every component.
Recall that we denote this space by $V^{C^3}(M)$. Thus a small perturbation of $f$ can be denoted by
$f_t= f + t g, t \in J=(-\epsilon, \epsilon)$, $ g \in V^{C^3}(M)$.

Let $\CB_m, m=0,1,2$ be $C^m(\mathbb T^n)$, the space consisting of functions with $i$th-order continuous partial derivatives.

The transfer  operator $\CL_t$ on $\CB_m$  is defined by
\begin{equation}   \label{transferopB}
\CL_t \varphi(x) = \sum_{y \in f_t^{-1}(x)}  \frac{\varphi(y)}{Jf_t(y)},
             \end{equation}
where $Jf_t$ is the Jacobian of $f_t$: $Jf_t =|\det(Df_t)|$.
Since the manifold $M=\T^n$ is compact, we see that there exists a small $\epsilon >0$ such that $\CL_t$ is a bounded linear  operator on $\CB_m$ for any $t\in (-\epsilon, \epsilon)$ and $i=0,1,2$.
{
Note that $f_t=f+t g$ is $C^\infty$ smooth in $t$.
Using inverse function theorem,
shrinking $\epsilon$ if necessary, one may assume that all $\{f_t\}_{t\in (-\epsilon, \epsilon)}$ have
the same number of inverse branches, say  $y^i_t(x)$ for $i=1, 2, \dots, \ell$,
such that each $y^i_t$ is also $C^\infty$ smooth in $t$.
Now we rewrite \eqref{transferopB} as
\begin{equation}
\label{transferopB branch}
\CL_t \varphi(x) = \sum_{i=1}^\ell  \varphi(y^i_t(x)) Jy^i_t(x).
\end{equation}
From this expression, it is easy to see that
\begin{itemize}
\item
The map $  t\mapsto \CL_t \in L(\CB_m, \CB_m)$ is $C^\infty$ smooth in $(-\epsilon, \epsilon)$;
\item
By taking derivative with respect to $t$ on both sides of \eqref{transferopB branch}, we get
$ (-\epsilon, \epsilon) \ni t\mapsto \p_t\CL_t \in L(\CB_m, \CB_{m-1})$
is also $C^\infty$ smooth.
\end{itemize}
It follows that for any $\eps_1\in (0, \eps)$,
the map
$ [-\epsilon_1, \epsilon_1] \ni t\mapsto \p_t\CL_t \in L(\CB_m, \CB_{m-1})$
is Lipschitz continuous.
}

Alternatively, $\CL_t$: $\CB_m \to \CB_m$ can be characterized by the following duality
\begin{equation}\label{duality}
\left< \CL_t\varphi, \psi\right> = \left< \varphi, \psi\circ f_t \right>,
\end{equation}
where $\left<\varphi, \psi\right>$ denote the inner product $\int_{\mathbb T^n} \varphi \psi d x $ and $dx$ is just the Lebesgue measure on $\mathbb T^n$.

Let $\rho_t$ denote the density function of the unique SRB measure of $f_t$. It is the eigenfunction of $\CL_t$ corresponding to the simple isolated eigenvalue $1$:  $\CL_t \rho_t =\rho_t$.

We now verify that the transfer operator $\CL_t$ satisfies  Assumptions (H1)-(H3) in previous subsections.

(H1)  By either the Lasota-York estimates or the invariant cone technique (see Baladi \cite{Ba18} Chapter 2, also Ruelle \cite{Ru89}, Liverani \cite{Li03}), one can show that acting on $\CB_m, m=1,2,$ there is some $\eta \in (0,1)$ such that for each $t$, the transfer operator $\CL_t$ has 1 as an isolated eigenvalue and a spectral gap of $1-\eta$.  By Theorem \ref{decomposition2} and the mixing property of $f_t$, we can write
$\CL_t=\CP_t+\CN_t$ such that

(1) $\CP_t$ is a one-dimensional projection defined by
\begin{equation}\label{Proj}
\CP_t\varphi =\rho_t \int \varphi d x,
\end{equation}
where $\rho_t$ is the SRB density of $f_t=f + t g$.

(2) $\CN_t$ is a bounded linear operator with spectral radius uniformly less than $r$.

(3) $\CP_t\CN_t =\CN_t\CP_t =0.$

\begin{comment}
No reference  gives the exact statements we need here. Ruelle's results are more general and proved in a different way. Baladi's is for dimension one map.
\end{comment}

(H2)  The differentiability of $\CL_t: \CB_m\to \CB_{m-1}, m=1,2$ with respect to $t$ can be seen from the definition of the transfer operator since the inverse branches of $f_t$ depend  smoothly on the parameter $t$. Nevertheless, the duality characterization is more convenient to derive the explicit formula of $\partial_t \CL_t$. Taking derivative with respect to $t$ on both sides of $\left< \CL_t\varphi, \psi\right> = \left< \varphi, \psi\circ f_t \right>$,
\[\left< (\partial_t \CL_t)\varphi, \psi\right>= \left<\varphi, \partial_t(\psi \circ f_t )\right> = \left<\varphi, \partial_t(\psi \circ (f+ t g) )\right>   \]
\[ = \left<\varphi, D\psi (f_t) \cdot g \right>= \left<\varphi g, D\psi \circ f_t\right>,\]
where $ \left<\varphi g, D\psi \circ f_t\right>  $ is the gradient of $\psi$ in the direction of the vector $\varphi g$ evaluated at $f_t$. Without loss of generality, we write $D\psi $ as a vector $(\frac{\partial \psi}{\partial x_1}, \frac{\partial \psi}{\partial x_2}, \cdots, \frac{\partial \psi}{\partial x_n} )$ and the vector field $g$ as $(g_1, g_2, \cdots, g_n).$
Thus, \[\left<\varphi g, D\psi \circ f_t\right>  = \sum_{k=1}^n    \left<    \varphi g_k,  \frac{\partial \psi}{\partial x_k} \circ f_t  \right>= \sum_{k=1}^n  \left< \CL_t   (\varphi g_k) , \frac{\partial \psi}{\partial x_k}  \right> = \left< \CL_t \varphi g, \grad \psi \right>.\]

Using the integration by parts over the closed manifold $M$, we have
\[ < \CL_t (\varphi g), \grad \psi > = - < \Div (\CL_t (\varphi g)), \psi> .\]
Thus, for any $\varphi, \psi  \in \CB_2$, we have
\[ \left< (\partial_t \CL_t)\varphi, \psi\right> =-  \left<   \Div (\CL_t (\varphi g)), \psi\right>,\]
where $\CL_t(\varphi g)$ is understood as $(  \CL_t \varphi g_1,  \CL_t \varphi g_2, \cdots, \CL_t \varphi g_n                    )$ when $g=(g_1, g_2, \cdots, g_n)$.
It yields that $\partial_t\CL_t$: $\CB_2 \to \CB_1$ (or, $\CB_1 \to \CB_0$) is given
by
\begin{equation}\label{dev.L}
( \partial_t\CL_t)\varphi = - \Div[\CL_t (\varphi g)].\end{equation}

(H3)  It is well known that  $\CL_t$ satisfies a uniform Lasota-Yorke estimates on the pair $(\CB_1, \CB_0)$. (H3) follows by applying the perturbation theory of Gouezel-Liverani (\cite{GL}, Section 8).

{
We notice that
the family of transfer operators $\CL_t$  for expanding endomorphisms
acting on the Banach spaces $\CB_m=C^m(\mathbb T^n), m=0,1,2$,
satisfies
the Lasota-Yorke type inequalities (A.2)-(A.7) in Appendix A.3 of Baladi's book \cite{Ba18}.
See \cite{Ba, Ba18} for more details about these Lasota-Yorke estimates.
Then (H3) follows from  the Gouezel-Keller-Liverani  perturbation theory,
as we have explained in the last paragraph of Subsection~\ref{subsec: assump}.
}

 \bigskip

Recall that $\rho_t$ is the density of the SRB measure of $f_t$. For convenience, we write $\CL=\CL_0$ and $\rho=\rho_0$. We now derive an explicit formula for the linear response function, i.e., the derivative formula of $\rho_t$ with respect to $t$. By  (\ref{Proj}),
we have $\rho_t=\CP_t\varphi$ for any $\varphi$ with $\int \varphi dx =1$, and hence, $\partial_t \rho_t=(\partial_t\CP_t)\varphi$.
In particular, if we let $\varphi= \rho_t=\CP_t\rho_t \in \CB_2$, then by Lemma \ref{PPP}, we have
\[  \partial_t\rho_t = (\partial_t\CP_t)\CP_t\rho_t=\sum_{n=0}^\infty \CL^n_t(\partial_t\CL_t)\CP^2_t\rho_t
= \sum_{n=0}^\infty \CL^n_t(\partial_t\CL_t) \rho_t .        \]
Therefore, we have the formula of the linear response function (in the direction of a vector field $g$)
\[\xi= \partial_t\rho_t |_{t=0}= \sum_{n=0}^\infty \CL^n(\partial_t\CL_t)|_{t=0} \rho = - \sum_{n=0}^\infty \CL^n \Div[\CL(g\rho)] \in \CB_0.\]

Finally, we have that for any smooth function $\psi$ on $M$, the function $t \to \int  \phi \rho_t dx $ is a differentiable function for any given $g \in V^3(M)$ and \[ \frac{d}{dt}\big|_{t=0} \int  \psi \rho_t dx  = - \sum_{n=0}^\infty \int  \psi   \CL^n \Div[\CL(g\rho)]  dx =-  \sum_{n=0}^\infty \int     \Div[\CL(g\rho)]\cdot  \psi \circ f^n dx. \]

{

\section{Perturbation Theory for Two-parameter Family of Operators }\label{sec.abs2}

\subsection{A two-parameter family of operators with a spectral gap on nested Banach spaces}

Consider a two parameter family of bounded linear operators $\CL_{\bt}$,
with $\bt=(t_1, t_2)\in \bJ:=J_1\times J_2\subset \R^2$,
compatibly acting on three nested Banach spaces $\CB_3\subset \CB_2 \subset \CB_1\subset \CB_0$,
such that the following assumptions (H1')-(H3') hold:

(H1'): There is $\eta< 1$ such that for all $\bt\in \bJ$,
$\CL_{\bt}|_{\CB_m}, m=1, 2, 3$ has a  $(1, \eta)$-spectral gap.

We denote the corresponding spectral gap decomposition $\CL_{\bt}=\CP_{\bt} + \CN_{\bt}$.

\medskip

(H2'):  Write $\bs=(s_1, s_2)\in \bJ$.
For any $\bt \in \bJ$, the limits
\beqn
\partial_i \CL_{\bt} := \lim_{s_i \to t_i } \frac{1}{s_i-t_i} ( \CL_{\bs} - \CL_{\bt})
\ \ \text{for} \ i=1, 2
\eeqn
 exists as a bounded linear operator from $\CB_k$ to $\CB_{k-1}$ for $k=2, 3$;
and the limit
\beqn
\partial_i \partial_j \CL_{\bt}:=
\lim_{s_i \to t_i } \frac{1}{s_i-t_i} ( \p_j \CL_{\bs} - \p_j\CL_{\bt})
\ \ \text{for} \ i, j=1, 2
\eeqn
 exists as a bounded linear operator from $\CB_3$ to $\CB_1$.

\medskip

Assume $\CL_{\bt}, \bt \in \bJ$, satisfies (H1') and let $\gamma$ be the circle centered at $1$ in the complex plane with a radius $\frac{1-\eta}{2}$ and $R_{\bt}(z)=\left(z-\CL_{\bt} \right)^{-1}$ be the resolvent. Assumption (H3') is about the continuity of the resolvent as an operator from $\CB_m$ to $\CB_{m-1}$ for $m=1, 2$.

(H3'):  For any $\bt \in\bJ$ and any $m=1, 2$,
\[\lim_{\bs \to \bt} \sup_{z \in \gamma} \|R_{\bs}(z) - R_{\bt}(z)\|_{L(\CB_m, \CB_{m-1})}=0.\]

\medskip

Applying similar computation in Appendix~\ref{sec.abs},
we obtain the following results.

\begin{lemma}
\label{merom'}
For any $\bt \in \bJ$ and any $z \in V\setminus \{1\}$,
where  $V=\{ z\in \mathbb C: |z|>\eta.\}$,
\[ R_{\bt}(z) = \frac{\CP_\bt}{z-1} + \CQ_\bt(z),\]
where
\begin{equation}\label{Qt'}
 \CQ_\bt(z) := R_\bt(z)(I - \CP_\bt) = (z- \CN_\bt)^{-1}(I-\CP_\bt)
 \end{equation}
 is a homomorphic function on $V$.
\end{lemma}

Note that for any $z\in V\setminus \{1\}$ with $|z-1|<1-\eta$,
\beqn
 (z- \CN_\bt)^{-1} =\left[ (z-1) + \left( I- \CN_\bt\right)\right]^{-1}
 =\sum_{n=0}^\infty (-1)^n (z-1)^n \left( I- \CN_\bt\right)^{-n-1},
\eeqn
and thus
\beqn
\CQ_{\bt}(z) =\sum_{n=0}^\infty \CQ_\bt^n (z-1)^n,
\ \ \text{where} \
\CQ_\bt^n:= (-1)^n \left( I- \CN_\bt\right)^{-n-1}(I-\CP_{\bt}).
\eeqn
Note that $\CP_\bt\CQ_\bt^n=\CQ_\bt^n\CP_\bt=0$
for any $n\ge 0$.
In particular,
\beq
\label{def Q01}
\CQ_{\bt}^0=(I-\CN_{\bt})^{-1} (I-\CP_{\bt})
\ \ \text{and} \ \
\CQ_{\bt}^1 =- (I-\CN_{\bt})^{-2}  (I-\CP_{\bt}).
\eeq

\medskip

\begin{lemma}\label{dev.R'}
For any $\bt \in \bJ$ and any $z \in \gamma$, we have
\beqn
\partial_i R_{\bt}(z):= \lim_{s_i \to t_i} \frac{1}{s_i-t_i}(R_{\bs}(z)-R_{\bt}(z) )
=R_{\bt}(z)(\partial_{i} \CL_{\bt} )R_{\bt}(z) , \
\ \text{for} \ i=1, 2
\eeqn
which is a bounded linear operator from $\CB_k$ to $\CB_{k-2}$ for $k=2, 3$;
and thus for $i, j=1, 2$,
\beqn
\begin{split}
\p_i \p_j R_{\bt}(z):=
&\lim_{s_i \to t_i} \frac{1}{s_i-t_i}(\p_j R_{\bs}(z)- \p_j R_{\bt}(z) ) \\
= & \, \, R_{\bt}(z)(\partial_i \CL_{\bt}) R_{\bt}(z)  (\partial_{j} \CL_{\bt} )R_{\bt}(z) \\
& + R_{\bt}(z)(\partial_j \CL_{\bt}) R_{\bt}(z)  (\partial_{i} \CL_{\bt} )R_{\bt}(z) \\
& +  R_{\bt}(z)(\p_i \partial_j \CL_{\bt}) R_{\bt}(z)
\end{split}
\eeqn
is a bounded linear operator from $\CB_3$ to $\CB_0$.
\end{lemma}

\begin{theorem}\label{dev.P'}
For any $ \bt \in\bJ$, we have
\beqn
 \partial_i\CP_{\bt} = \CP_{\bt}(\partial_i \CL_{\bt})\CQ_{\bt}^0  + \CQ_{\bt}^0(\partial_i \CL_{\bt}) \CP_{\bt},
 \ \text{for}\ i=1, 2,
\eeqn
which is  bounded linear operator from $\CB_k$ to $\CB_{k-2}$ for $k=2, 3$;
also, for $i, j=1, 2$,
\beqn
\begin{split}
\p_i \partial_j\CP_{\bt}
&  =\CP_{\bt}(\partial_i \CL_{\bt})\CQ_{\bt}^0  (\partial_j\CL_{\bt})\CQ_{\bt}^0
+\CQ_{\bt}^0 (\partial_i \CL_{\bt})\CP_{\bt}  (\partial_j\CL_{\bt})\CQ_{\bt}^0
+\CQ_{\bt}^0(\partial_i \CL_{\bt})\CQ_{\bt}^0  (\partial_j\CL_{\bt}) \CP_{\bt}
\\
& + \CP_{\bt}(\partial_i \CL_{\bt})\CP_{\bt}^0  (\partial_j\CL_{\bt})\CQ_{\bt}^1
+\CP_{\bt} (\partial_i \CL_{\bt})\CQ_{\bt}^1  (\partial_j\CL_{\bt})\CP_{\bt}
+\CQ_{\bt}^1 (\partial_i \CL_{\bt})\CP_{\bt}  (\partial_j\CL_{\bt}) \CP_{\bt}
\\
&  +\CP_{\bt}(\partial_j \CL_{\bt})\CQ_{\bt}^0  (\partial_i\CL_{\bt})\CQ_{\bt}^0
+\CQ_{\bt}^0 (\partial_j \CL_{\bt})\CP_{\bt}  (\partial_i\CL_{\bt})\CQ_{\bt}^0
+\CQ_{\bt}^0(\partial_j \CL_{\bt})\CQ_{\bt}^0  (\partial_i\CL_{\bt}) \CP_{\bt}
\\
& + \CP_{\bt}(\partial_j \CL_{\bt})\CP_{\bt}^0  (\partial_i\CL_{\bt})\CQ_{\bt}^1
+\CP_{\bt} (\partial_j \CL_{\bt})\CQ_{\bt}^1  (\partial_i\CL_{\bt})\CP_{\bt}
+\CQ_{\bt}^1 (\partial_j \CL_{\bt})\CP_{\bt}  (\partial_i\CL_{\bt}) \CP_{\bt}
\\
& + \CP_{\bt}(\partial_i \p_j \CL_{\bt})\CQ_{\bt}^0  + \CQ_{\bt}^0(\partial_i \p_j \CL_{\bt}) \CP_{\bt}
\end{split}
\eeqn
which is  bounded linear operator from $\CB_3$ to $\CB_{0}$.
\end{theorem}

We shall focus on the composition $(\partial_i \CP_\bt) \CP_\bt$
and $(\partial_i \p_j\CP_\bt) \CP_\bt$.

\begin{lemma}
\label{lem: p_ij formula}
For any $\bt \in \bJ$, we have
\beqn
(\partial_i \CP_\bt) \CP_\bt =
\CQ_{\bt}^0(\partial_i \CL_{\bt}) \CP_{\bt}
\ \ \text{for}\ i=1, 2.
\eeqn
Also, for $i, j=1, 2$, we have
\beqn
\begin{split}
(\partial_i \p_j \CP_\bt) \CP_\bt
& =\CQ_{\bt}^0(\partial_i \CL_{\bt})\CQ_{\bt}^0  (\partial_j\CL_{\bt}) \CP_{\bt}
+\CP_{\bt} (\partial_i \CL_{\bt})\CQ_{\bt}^1  (\partial_j\CL_{\bt})\CP_{\bt}
+\CQ_{\bt}^1 (\partial_i \CL_{\bt})\CP_{\bt}  (\partial_j\CL_{\bt}) \CP_{\bt}
\\
& +\CQ_{\bt}^0(\partial_j \CL_{\bt})\CQ_{\bt}^0  (\partial_i\CL_{\bt}) \CP_{\bt}
+\CP_{\bt} (\partial_j \CL_{\bt})\CQ_{\bt}^1  (\partial_i \CL_{\bt})\CP_{\bt}
+\CQ_{\bt}^1 (\partial_j \CL_{\bt})\CP_{\bt}  (\partial_i \CL_{\bt}) \CP_{\bt}
\\
& + \CQ_{\bt}^0(\partial_i \p_j \CL_{\bt}) \CP_{\bt}.
\\
\end{split}
\eeqn
\end{lemma}

\subsection{Application to Toral Expanding Endomorphisms}\label{sec.abs3}

Recall that $E^{C^r}(M)$ is the space of expanding endomorphisms of a closed Riemannian manifold $M$.
For simplicity, we only consider the $n$-torus $M=\T^n\cong \R^n/\mathbb{Z}^n$,
and we assume that $r\ge 4$ (in the previous sections, we only assume $r\ge 3$).
Note that $E^{C^r}(M)$ is a Banach manifold whose tangent space
can be identified as the space $V^{C^r}(M)$
of maps from $\R^n$ to itself that are $1$-periodic in every component.

Recall that the transfer operator of $f\in  E^{C^r}(M)$ is denoted by $\CL_f$,
which is compatibly defined as bounded linear operators
of the nested Banach spaces $\CB_k=C^k(M)$ for $k=0,1,2, 3$.
Using Kato's perturbation theory  and locally uniform Lasota-Yorke estimates,
we can show that the spectral gap decomposition is locally uniform.

\begin{lemma}
\label{lem: sp gap uniform}
Let $f\in  E^{C^r}(M)$. There are
$\eps_f>0$,  $\eta_f\in (0, 1)$
and  $K_f>0$ such that if $h\in E^{C^r}(M)$ satisfies $\|h-f\|_{C^r}<\eps_f$,
then its transfer operator $\CL_h$ has a $(1, \eta_f)$-spectral gap.
Moreover,
denote the corresponding spectral gap decomposition by $\CL_h=\CP_h+\CN_h$,
then for any $m=1, 2, 3$ and any $n\in \mathbb{N}$,
\beqn
\|\CL_h\|_{\CB_m\to \CB_m}\le K_f, \
\|\CP_h\|_{\CB_m\to \CB_m}\le K_f, \
\text{and} \
\|\CN_h^n\|_{\CB_m\to \CB_m} \le K_f \eta_f^n.
\eeqn
\end{lemma}

We are particularly interested in
a  two-parameter family  of the following form:
\beq
\label{def f_bt}
f_{\bt}=f+ g_0 + t_1 g_1 + t_2 g_2,
\eeq
where $f\in  E^{C^r}(M)$,
$g_0, g_1, g_2\in E^{C^r}(M)$ are a priori given, and
$\bt=(t_1, t_2) \in \bJ=J_1\times J_2 \subset \R^2$
is regarded as two-dimensional parameter vector.
We assume that
all $f_\bt$ lies in the $C^r$ neighborhood of $f$ of size $\eps_f$,
where  $\eps_f$ is given by Lemma~\ref{lem: sp gap uniform}.
Denote the transfer operator of $f_\bt$ by $\CL_\bt$,
with the
spectral gap decomposition by $\CL_\bt=\CP_\bt+\CN_\bt$.
It follows from Lemma~\ref{lem: sp gap uniform} that
 for any $m=1, 2, 3$, any $n\in \mathbb{N}$ and any $\bt\in \bJ$ that
\beq
\label{eq: uniform sp gap}
\|\CL_\bt\|_{\CB_m\to \CB_m}\le K_f, \
\|\CP_\bt\|_{\CB_m\to \CB_m}\le K_f, \
\text{and} \
\|\CN_\bt^n\|_{\CB_m\to \CB_m} \le K_f \eta_f^n.
\eeq
Recall the definitions of
$\CQ_{\bt}^0$ and $\CQ_{\bt}^1$ are given by \eqref{def Q01},
then \eqref{eq: uniform sp gap} immediately implies that
there is $K_f'>0$ such that for all $\bt\in \bJ$,
\beq
\label{eq: uniform sp gap1}
\|\CQ_{\bt}^0\|_{\CB_m\to \CB_m}\le K_f',
\ \text{and} \
\|\CQ_{\bt}^1\|_{\CB_m\to \CB_m}\le K_f',
\ \text{for} \ m=1, 2, 3.
\eeq

Applying similar computation in Subsection~\ref{app: linear response},
we have precise formula for the partial derivatives
$\p_i \CL_\bt$ and $\p_i\p_j \CL_\bt$.
For the vector valued functions
$g_1, g_2\in E^{C^r}(M)$, we write
$g_1=(g_1^1, \dots, g_1^n)$
and
$g_2=(g_2^1, \dots, g_2^n)$.

\begin{lemma}
\label{lem: par der CL_bt}
Let $\CL_\bt$ be the family of transfer operator corresponding to
$f_\bt$ given in \eqref{def f_bt}.
Then for any $\varphi\in \CB_m$ with $m=1, 2, 3$,
\beqn
\p_i \CL_\bt \varphi = - \sum_{k=1}^n \frac{\p}{\p x_k} \left[ \CL_\bt \left( \varphi g_i^k\right)\right]
= -\Div  \left[ \CL_\bt \left( \varphi g_i\right)\right].
\eeqn
Also, for any $\varphi\in \CB_m$ with $m=1, 2$,
\beqn
\p_i  \p_j \CL_\bt \varphi
= \sum_{k=1}^n \sum_{\ell=1}^n \frac{\p^2}{\p x_k \p x_\ell} \left[ \CL_\bt \left( \varphi g_i^k g_j^\ell\right)\right].
\eeqn
\end{lemma}

It follows from Lemma~\ref{lem: par der CL_bt} and Inequality~\eqref{eq: uniform sp gap} that for any
$\varphi\in \CB_m$ with $m=1, 2, 3$,
\beqn
\|\p_i \CL_\bt \varphi \|_{\CB_{m-1}}
\le n \max_{1\le k\le n} \left\| \CL_\bt \left( \varphi g_i^k\right)  \right\|_{\CB_{m}}
\le n K_f \max_{1\le k\le n} \left\|  \varphi g_i^k  \right\|_{\CB_{m}}
\le n K_f \left\|  \varphi  \right\|_{\CB_{m}}\left\|   g_i   \right\|_{\CB_{m}},
\eeqn
that is,
\beq
\label{eq: est p_i CL}
\|\p_i \CL_\bt  \|_{\CB_m\to \CB_{m-1}} \le n K_f \left\|   g_i   \right\|_{\CB_{m}}.
\eeq
Similarly, for $m=2, 3$,
\beq
\label{eq: est p_ij CL}
\|\p_i\p_j \CL_\bt  \|_{\CB_m\to \CB_{m-2}} \le n^2 K_f \left\| g_i \right\|_{\CB_{m}}\left\| g_j \right\|_{\CB_{m}}.
\eeq

\medskip

Let $\rho_\bt$ be the  SRB density  of $f_\bt$, then it has the representation
$\rho_\bt=\CP_\bt \varphi$ for any $\varphi\in \CB_3$ with $\int \varphi dx=1$.
In particular, $\rho_\bt=\CP_\bt 1$, where $1$ is the constant one function,
then by Inequality~\eqref{eq: uniform sp gap}, we have
\beqn
\|\rho_\bt\|_{C^0}\le \|\CP_\bt 1\|_{\CB_3}
\le \|\CP_\bt\|_{\CB_3\to \CB_3} \|1\|_{\CB_3}\le K_f.
\eeqn
On the other hand, for any $i=1, 2$,
we have $\p_i \rho_\bt = \p_i\CP_\bt \varphi$, and by taking $\varphi=\CP_t 1$, we get
$\p_i \rho_\bt=(\p_i\CP_\bt)  \CP_\bt 1$. By Lemma~\ref{lem: p_ij formula}
and Inequalities~\eqref{eq: uniform sp gap}\eqref{eq: uniform sp gap1}\eqref{eq: est p_i CL}, we have
\beqn
\|\p_i \rho_\bt\|_{C^0}
\le \|(\p_i\CP_\bt)  \CP_\bt 1\|_{\CB_2}
\le \|\CQ_{\bt}^0 \|_{\CB_2\to \CB_2}
\| \partial_i \CL_{\bt}\|_{\CB_3\to \CB_2} \| \CP_{\bt}\|_{\CB_3\to \CB_3} \|1\|_{\CB_3}
\le n K_f^2 K_f' \left\|   g_i   \right\|_{\CB_3}.
\eeqn
Similarly, By the second identity of Lemma~\ref{lem: p_ij formula}
and Inequalities~\eqref{eq: uniform sp gap}\eqref{eq: uniform sp gap1}\eqref{eq: est p_i CL}~\eqref{eq: est p_ij CL}, we have
\beqn
\|\p_i \p_j \rho_\bt\|_{C^0}
\le \|(\p_i \p_j \CP_\bt)  \CP_\bt 1\|_{\CB_1}
\le 6n^2 K_f^3 (K_f')^2 \left\|g_i \right\|_{\CB_3}\left\|g_j \right\|_{\CB_3}
+ n^2 K_f^2 K_f'\left\|g_i \right\|_{\CB_3}\left\|g_j \right\|_{\CB_3}.
\eeqn

We summarize by taking
$
C_f=\max\{K_f, n K_f^2 K_f',  6n^2 K_f^3 (K_f')^2+ n^2 K_f^2 K_f'\}
$.

\begin{theorem}
\label{thm: der rho bt}
Let $f_\bt$ be the family of expanding endomorphisms given in \eqref{def f_bt},
and $\rho_\bt$ be the corresponding SRB density  of $f_\bt$.
Then there is a constant $C_f>0$ such that for all $\bt\in \bJ$ and any $i, j=1, 2$,
\beqn
\|\rho_\bt\|_{C^0}\le C_f, \ \
\|\p_i \rho_\bt\|_{C^0} \le C_f\left\|g_i\right\|_{C^r},  \ \
\|\p_i \p_j \rho_\bt\|_{C^0} \le C_f \left\|g_i\right\|_{C^r}\left\|g_j\right\|_{C^r}.
\eeqn
\end{theorem}

}

\begin{comment}
Given $\alpha\ge 1$,
we denote by $\EEnd^\alpha(\T^d)$
the space of all $C^\alpha$ smooth uniformly expanding endomorphisms of
the $d$-dimensional torus.
Let $\CL$ be the transfer operator of $f$,
then it is well known (see e.g.~\cite{Ba}) that the ACIP density function $\rho$ is
the unique fixed point of   $\CL$.
Moreover,
For any $\beta\in (0, \alpha-1]$,
the action of $\CL$ restricted on $C^{\beta}(\T^d)$ has \emph{spectral gap},
i.e.,
$\mathrm{Spectrum}(\CL)$ consists of a simple leading eigenvalue $1$
and a compact subset  strictly inside the unit disc.
Equivalently,
we can write
$\CL=\CP + \CN$
with the following properties:
\begin{enumerate}
\item $\CP$ is a one-dimensional projection, i.e., $\CP^2=\CP$ and $\dim(\mathrm{Im}(\CP))=1$;
\item $\CN$ is a bounded linear operator with spectral radius strictly less than $1$;
\item $\CP\CN=\CN\CP=0$.
\end{enumerate}
Furthermore, $\CP$ is the leading eigen-projection of $\CL$ given by the following formula:
\beq
\label{def CP0}
\CP\varphi =\rho \int \varphi \, dm,
\ \ \text{for any} \ \
 \varphi\in C^{\beta}(\T^d).
\eeq

\medskip

Q1:  Examples of different norms on $\CB_m$:  $C^r$ norm, $L^p$ norm, Sobolev norm, etc. They lead to different statements on differentiability of SRB.

Q2: Proofs of Lasota-Yorke in different norms.

\end{comment}


\begin{thebibliography}{99}


\bibitem[AF]{AF} R.A. Adama and J.J.F. Fournier, Sobolev Spaces, 2nd ed. 2003, Elsevier Science Ltd.

\bibitem[Ba]{Ba} V. Baladi, Positive Transfer Operators and Decay of Correlations, World Scientific, Singapore, New Jersey, London Hong Kong, 2000

 \bibitem[Ba18]{Ba18} V.   Baladi, Viviane,
Dynamical zeta functions and dynamical determinants for hyperbolic maps.
A functional approach
Ergeb. Math. Grenzgeb. (3), 68 [Results in Mathematics and Related Areas. 3rd Series. A Series of Modern Surveys in Mathematics]
Springer, Cham, 2018. xv + 291 pp.
	
\bibitem[Belse]{Belse} V. Baladi, Linear response or else. In: (English summary) Proceedings of the International Congress of
Mathematicians-Seoul, vol. III, pp. 525–545, Kyung Moon Sa, Seoul (2014)



\bibitem[G96]{G96} G. Gallavotti,   Chaotic hypothesis: Onsager reciprocity and fluctuation-dissipation theorem.
J. Statist. Phys. 84 (1996), no. 5-6, 899-925

\bibitem[G06]{G06}G. Gallavotti,   Entropy, thermostats, and chaotic hypothesis. Chaos 16 (2006), no. 4, 043114, 6 pp

\bibitem[GC]{GC} G. Gallavotti and E.G.D.  Cohen,  Dynamical ensembles in stationary states. J. Statist. Phys. 80 (1995), no. 5- 6, 931-970.

\bibitem[Gi]{Gi} S. Giuntini  A Remark on Modified Euler's Method for Differential Equations in Banach Spaces.
Universitatis Iagellonicae ACTA Mathematica 1985.

\bibitem[GL]{GL} S. Gou\"ezel and C. Liverani,  Banach spaces adapted to
Anosov systems. Ergodic Theory Dynam. Systems, 26(1):189-217, 2006

\bibitem[He]{He} E. Hebey,    Nonlinear Analysis on Manifolds: Sobolev spaces and inequalities, AMS,       1999

 \begin{comment}
\bibitem[Jaynes]{Ja}  Jaynes, MEPP (to be completed)  E.T. Jaynes, “Information Theory and Statistical Mechanics”, Physical Review, Vol. 106, No. 4, 1957, pp. 620–630

Beyond the Second Law

Note: Jaynes did not propose the Maximal Entropy Production Principle. His work on Maximal Entropy Principle inspired the later formulation of MEPP  for  open systems.  For evolution along the gradient of the entropy see next reference

\end{comment}

\bibitem[JP]{JP}  A. Jane$\check{\rm c}$ka  and M. Pavelka,   Gradient Dynamics and Entropy Production Maximization,  Journal of Non-Equilibrium Thermodynamics, vol. 43, no. 1, 2018, pp. 1-19. https://doi.org/10.1515/jnet-2017-0005




\bibitem[J12]{J12} M. Jiang,    Differentiating potential functions of SRB measures on hyperbolic attractors, Ergodic Theory Dynam. Systems 32(4)  (2012) , 1350 - 1369

\bibitem[J21]{J21} M. Jiang,    Chaotic hypothesis and the second law of thermodynamics. Pure Appl. Funct. Anal. 6 (2021), no. 1, 205 - 219


\bibitem[JL22]{JL22} M. Jiang, SRB entropy of Markov Transformations, J. Stat. Physics, 188 (2022) No.3 Paper No. 24, 18 pp

\bibitem[J24]{J24} M. Jiang,   Gradient flow of the Sinai-Ruelle-Bowen entropy,
 Comm. Math. Phys. 405 (2024), no. 5, Paper No. 118, 22 pp.

\bibitem[LY]{LY} A. Lasota and J.A.  Yorke,  The generic property of existence of solutions of differential equations in Banach space. J. Differential Equations 13 (1973), 1-12.

\bibitem[Li03]{Li03} C. Liverani, Carlangelo,
Invariant measures and their properties. A functional analytic point of view. (English summary) Dynamical systems. Part II, 185–237.
Pubbl. Cent. Ric. Mat. Ennio Giorgi [Publications of the Ennio de Giorgi Mathematical Research Center]
Scuola Normale Superiore, Pisa, 2003

\bibitem[Maas]{Maas} J. Maas, Gradient flows of the entropy for finite Markov chains, J. Func. Analysis 261 (2011) 2250-2292

\bibitem[Kato]{Kato}   T. Kato. Perturbation theory for linear operators. Classics in Mathematics.
Springer-Verlag, Berlin, 1995. Reprint of the 1980 edition

\bibitem[KH]{KH} A. Katok and B.  Hasselblatt,
Introduction to the modern theory of dynamical systems.
With a supplementary chapter by Katok and Leonardo Mendoza. Encyclopedia of Mathematics and its Applications, 54. Cambridge University Press, Cambridge, 1995

\bibitem[KL]{KL} G. Keller and C. Liverani, Stability of spectrum for transfer operators, Ann. Scuola Norm. Sup. Pisa Cl. Sci. 28, 141-152, 1999

\bibitem[Mane]{Mane} R. Ma\~n\'e,   Ergodic theory and differentiable dynamics. Translated from the Portuguese by Silvio Levy. Ergebnisse der Mathematik und ihrer Grenzgebiete (3) [Results in Mathematics and Related Areas (3)], 8. Springer-Verlag, Berlin, 1987.


\bibitem[Ru89]{Ru89}  D. Ruelle,   The thermodynamic formalism for expanding maps. Comm.
Math. Phys., 125(2):239-262, 1989

\bibitem[Ru97]{Ru97}  D. Ruelle,   Differentiation of SRB states. Comm. Math. Phys. 187 (1997),
no. 1, 227-241.


\bibitem[Ru03]{Ru03}  D. Ruelle,
 Correction and complements: Differentiation of SRB states,
Comm. Math. Phys. 234 (2003), no. 1, 185-190


\bibitem[Young]{Young} L-S. Young,  What are SRB measures, and which dynamical systems have them? Dedicated to David Ruelle and Yasha Sinai on the occasion of their 65th birthdays. J. Statist. Phys. 108 (2002), no. 5 -6, 733-754.

\end{thebibliography}
\end{document}